\documentclass{amsart}
\usepackage[margin=1.8cm]{geometry}
\usepackage{amssymb}
\usepackage{relsize}
\newcommand{\ccro}{{\cc^{\rho}}}

\newcommand{\ccpsro}{{\cc^{\ro}_s}}

\newcommand\numeq[1]%
  {\stackrel{\scriptscriptstyle(\mkern-1.5mu#1\mkern-1.5mu)}{=}}

\newcounter{relctr} %% <- counter for relations
\everydisplay\expandafter{\the\everydisplay\setcounter{relctr}{0}} %% <- reset every eq
\AtBeginDocument{} %% <- store original definition

\newcommand{\jo}{{j_1}}
\newcommand{\jtw}{{j_2}}

\usepackage{tikz-cd}
\usepackage{latexsym,amssymb,amsmath}
\usepackage{comment}
\usepackage{mathrsfs}
\usepackage{amsmath,amscd}
\usepackage[enableskew]{youngtab}
\usepackage[all]{xy}
\usepackage{empheq}

\usepackage{amsmath}

\usepackage{empheq}
\usepackage{xcolor}
\definecolor{lightgreen}{HTML}{90EE90}

\newcommand{\lag}{\langle}
\newcommand{\rag}{\rangle}
\newtheorem{theorem}{Theorem}[section]

\newtheorem{remark}[theorem]{Remark}
\newtheorem{notations}[theorem]{Notations}
\newcommand{\sent}{\mapsto}

\newcommand\C{\mathcal{C}}
\DeclareMathOperator{\id}{id}

\excludecomment{verlong}
\includecomment{vershort}
\excludecomment{noncompile}

\usepackage{tikz}
\usetikzlibrary{matrix}
\usepackage{combelow}
\usepackage{amsmath,amsthm,amsfonts,amssymb}
\usepackage{amsmath, amsthm, amssymb, amscd, amsfonts}
\usepackage[all]{xy}
\usepackage{amssymb, amsthm, amsmath}
\usepackage{amsfonts}
\usepackage{color}

\newcommand{\ra}{\rightarrow}

\newcommand{\ot}{\otimes}

\newcommand{\xra}{\xrightarrow}
\newcommand{\mtc}{\mathcal}

\newcommand{\lam}{\lambda}

\newcommand{\Lam}{\Lambda}

\newcommand{\al}{\alpha}
\newcommand{\eps}{\epsilon}

\newcommand{\ul}{\underline}

\newcommand{\lh}{\leftharpoonup}

\newtheorem{lem}[theorem]{Lemma}%[equation]
\theoremstyle{plain}

\newtheorem{thm}[theorem]{Theorem}%[equation]
\newtheorem{prop}[theorem]{Proposition}%[equation]
\newtheorem{defn}[theorem]{Definition}%[equation]
\newtheorem{cor}[theorem]{Corollary}
\newtheorem{rem}[theorem]{Remark}
\newcommand{\ch}{\chi}
\newcommand{\mtr}{\mathrm}

\newcommand{\ncm}{\newcommand}
\ncm{\np}{\newpage}
\ncm{\ebl}{\end{thebibliography}}
\ncm{\bbl}{\begin{thebibliography}}
\ncm{\chd}{_{ _{\ch}}}
\ncm{\ald}{_{ _{\al}}}
\newcommand{\blam}{\Lam}
\ncm{\cP}{\mathcal{P}}
\ncm{\ei}{e_i}
\ncm{\eij}{e_{i,\;j}}
\ncm{\bt}{\begin{thm}}
\ncm{\bdef}{\begin{defn}}
\ncm{\edf}{\end{defn}}
\ncm{\et}{\end{thm}}
\ncm{\bc}{\begin{cor}}
\ncm{\bl}{\begin{lem}}
\ncm{\el}{\end{lem}}
\ncm{\bpf}{\begin{proof}}
\ncm{\epf}{\end{proof}}
\ncm{\ec}{\end{cor}}
\ncm{\ord}{\mtr{ord}}
\ncm{\er}{\end{rem}}
\ncm{\br}{\begin{rem}}
\ncm{\bn}{\begin}

\ncm{\bp}{\begin{prop}}
\ncm{\ep}{\end{prop}}
\ncm{\bd}{
\begin{document}}
\ncm{\ed}{\end{document}}
\ncm{\beq}{\begin{equation}}
\ncm{\beqn}{\begin{equation*}}
\ncm{\eeq}{\end{equation}}
\ncm{\eeqn}{\end{equation*}}
\ncm{\bea}{\begin{eqnarray}}
\ncm{\eea}{\end{eqnarray}}
\ncm{\beanon}{\begin{eqnarray*}}
\ncm{\eeanon}{\end{eqnarray*}}\ncm{\ek}{\eps|_K}\ncm{\diez}{\#}
\ncm{\bwt}{\bowtie}
\ncm{\cC}{\mtc{C}}\ncm{\cc}{\mtc{C}}
\ncm{\cX}{\mtc{X}}
\ncm{\wt}{\widetilde}
\ncm{\sg}{\sigma}
\ncm{\Rep}{\mathrm{Rep}}
\DeclareMathOperator{\Irr}{Irr}
%\ncm{\Irr}{\mathrm{Irr}}
\ncm{\X}{\mathcal{X}}
\ncm{\cA}{\mathcal{A}}
\ncm{\HKer}{\mtr{HKer}}
\ncm{\LKER}{\mtr{LKer}}
\ncm{\aad}{\mtr{ad}}
\newcommand{\mbf}{\mathbb F}
\ncm{\Dr}{\mtr{D}}
\ncm{\cD}{{\mathcal{D}}}\ncm{\cd}{{\mathcal{D}}}\ncm{\ce}{{\mathcal{E}}}
\ncm{\G}{\mathcal{G}}
\ncm{\Dc}{\mtc{D}}
\ncm{\E}{\mtc{E}}
\ncm{\fp}{\mtr{FPdim}}
\ncm{\FPdim}{\mtr{FPdim}}
\ncm{\Vc}{\mtr{Vec}}
\ncm{\cK}{\mtc{K}}
\ncm{\cM}{\mtc{M}}
\ncm{\cE}{\mtc{E}}
\ncm{\cS}{\mtc{S}}

\newcommand{{\ipr}}{i'}
\newcommand{\tomega}{\widetilde{\omega}}

\DeclareMathOperator{\End}{End}
%\ncm{\End}{\mtr{End}}
%\ncm{\End}{\mtr{End}}
\ncm{\cop}{\mtr{cop}}
\ncm{\op}{\mtr{op}}
\ncm{\chr}{character }\ncm{\ck}{\mtc{K}}
\ncm{\bw}{\bwt}
%\ncm{\Vec}{\mtr{Vec}}
%\ncm{\cY}{\mtc{Y}}
\ncm{\hker}{\mtr{HKer}}
\ncm{\bx}{\boxtimes}%\ncm{\cd}{\cD}
\ncm{\blue}{\textcolor[rgb]{.00, .00, 1.00}}
\ncm{\red}{\textcolor[rgb]{1.00, .00, .00}}
\ncm{\green}{\textcolor[rgb]{.50, 0.20, .90}}
\ncm{\bne}{\begin{enumerate}}
\ncm{\ene}{\end{enumerate}}
\ncm{\lker}{\mtr{LKer}}
\ncm{\md}{\medbreak}
\ncm{\rep}{\Rep}\ncm{\ind}{\mtr{ind}}
\ncm{\mdn}{\md\noindent}
\ncm{\dd}{$}
\ncm{\up}{^}
\newcommand{\tcs}{\text}
\newcommand{\mbb}{\mathbb B}
\newcommand{\vs}{\mathbb V}%{\mathrm{Supp}}
\newcommand{\sth}{suppose that\;}
\newcommand\rad{\operatorname{rad}}
\newcommand{\itm}{\item}
\newcommand{\dbd}{$$}
%\linespread{0.90}
\newcommand{\mol}{\mtr{mod}}
 \newcommand{\ro}{\rho}
\newcommand{\irr}{\mathrm{Irr}}
\newcommand{\mbc}{\mathbb C}
\newcommand{\mbs}{\mathbb S}
\newcommand{\mbz}{\mathbb Z}
\newcommand{\ct}{\mtc T}%group action
\newcommand{\sm}{\setminus}
\newcommand{\epl}{^{+}}
\newcommand{\sbsq}{\subseteq}
\newcommand{\sbs}{\subset}
\newcommand{\cco}{\mtr{co}}
\newcommand{\cz}{\mathcal{Z}}
\newcommand{\dual}{^{*}}
\newcommand{\Gm}{\Gamma}
\ncm{\cY}{\mtc{Y}}
\newcommand\ZZ{{\mathbb Z}} 
\newcommand{\bab}{\color{DarkOrchid}{}}
\newcommand{\eab}{\normalcolor{}}
\newcommand{\subs}{\subsection}
\newcommand{\cv}{\mtc{V}}
%added april 20 2017%%%%%%%%
%%%%%%%%%%%%%%%%%%
%\newcommand{\tcs}{\text}
%\newcommand{\mbb}{\mathbb B}
%\newcommand{\mbs}{\mathbb S}
%\newcommand{\vs}{\mathbb V}
%{\mathrm{Supp}}
%\newcommand{\subs}{\subsection}
  \newcommand{\grn}{\green}
\newcommand{\dt}{\delta}

\newcommand{\ccf}{\mathrm{ {CF}(\cc)}}
\newcommand{\cce}{\mathrm{ {CE}(\cc)}}
\newcommand{\cecc}{\mathrm{ {CE}(\cc)}}
\newcommand{\cecd}{\mathrm{ {CE}(\cd)}}
\newcommand{\kk}{\Bbbk}
\newcommand{\otL}{\ot_{L}}
\newcommand{\otl}{\ot_{L}}
\newcommand{\unpsi}{1_{\psi}}
\newcommand{\epsi}{e_{\psi}}
\newcommand{\ephi}{e_{\phi}}
\newcommand{\ech}{e_{\ch}}
\newcommand{\nleftcid}{\text{left normal  coideal subalgebra}}
\newcommand{\dimL}{\dim_{\kk}L}
\newcommand{\cl}{\mtc L}
\newcommand{\mj}{\mtc J}
\newcommand{\tl}{\tilde L}
\newcommand{\tL}{\tilde L}
\newcommand{\tpsi}{\tilde(\psi)}
\newcommand{\tmx}{\tilde{\mtc X}}
\newcommand{\zlh}{\mathrm{ZL}}
\newcommand{\ba}{\mathrm A}
\newcommand{\bv}{\mathrm V}
\newcommand{\zhopf}{\mtc{Z}_{\mtr{Hopf}}}
\newcommand{\lstar}{L^{*}}
\newcommand{\ldstar}{L^{**}}
\newcommand{\mstar}{M^{*}}
\newcommand{\mdstar}{M^{**}}
\newcommand{\lkera}{\lker_{A}}
\newcommand{\mdprime}{M''}
\newcommand{\ldprime}{L''}
\newcommand{\cm}{\mtc M}
\newcommand{\ccm}{\mathcal M}
\newcommand{\cn}{\mathcal N}
\newcommand{\ccn}{\mathcal N}
\newcommand{\rx}{\mtr{Rex}}
\newcommand{\cca}{\ca}
\newcommand{\ih}{\underline{\mtr{Hom}}}
\newcommand{\cih}{\underline{\mtr{coHom}}}
\newcommand{\hm}{\mtr{ {Hom}}}
\newcommand{\cov}{\mtr{coev}}
\newcommand{\rora}{\rho^{\mtr{ra}}}
\newcommand{\rola}{\rho^{\mtr{la}}}
\newcommand{\cx}{\mtc X}
 \newcommand{\cZ}{\cz}
 \newcommand{\ca}{\cA}
 \newcommand{\stat}{\noindent}
 \newcommand{\bfa}{{\bf A}}
 \newcommand{\unu}{\mathbf{1}}
 \newcommand{\barzu}{{\bar {  Z}(\unu)}}
 
\newcommand{\idx}{\id_X}
\newcommand{\lprime}{L'}
\newcommand{\mprime}{M'}
\newcommand{\nat}{ \mtr{{  Nat}}}
\newcommand{\ft}{\mtc F_\lam}
\newcommand{\rhau}{\rightharpoonup}
\newcommand{\lhau}{\leftharpoonup}
\newcommand{\cf}{\mathrm{ {CF}}}

\newcommand{\cfc}{\mathrm{{CF}}(\cc)}
\newcommand{\csu}{\overline{\mathfrak{  C}}}
\newcommand{\cfcc}{\mathrm{ {CF}}(\cc)}%{\mathrm{ {Ch}}(\cc)}
\newcommand{\cfcd}{\mathrm{CF}(\cd)}
\newcommand{\cfd}{\mathrm{CF}(\cd)}
\newcommand{\czcc}{{\cz(\cc)}}
\newcommand{\czcd}{{\cz(\cd)}}
\newcommand{\czt}{{\cz(\cz(\cc))}}
\newcommand{\enx}{\mtr{  End}}
\newcommand{\runu}{R(\unu)}
%%%%ostrik-kirillov

\newcommand{\bdfn}{\bn{defn}}
\newcommand{\edfn}{\end{defn}}
%pivotal structure
\newcommand{\deltax}{\delta_X}
\newcommand{\deltav}{\delta_V}
%category of modules
\newcommand{\repcca}{\rep_\cc(A)}
%tensor products
\newcommand{\xotay}{X \ot_A Y}%tensprod over a
\newcommand{\xoty}{X \ot Y}
\newcommand{\votw}{V \ot W}
\newcommand{\votaw}{V \ot_A W}%tensprod over a
%dimensions
\newcommand{\dimax}{\dim_AX}
\newcommand{\dimccx}{\dim_\cc(X)}
\newcommand{\dimcca}{\dim_\cc(A)}
\newcommand{\dimccv}{\dim_\cc(V)}
\newcommand{\dima}{\dim_A}
\newcommand{\biga}{A}
\newcommand{\comp}{\mathbb C}
\newcommand{\tehtaa}{\theta_A}
\newcommand{\tetaa}{\theta_A}
\newcommand{\ida}{\id_A}
\newcommand{\hma}{\hm_A}
\newcommand{\hmcc}{\hm_\cc}
\newcommand{\fv}{F(V)}
\newcommand{\fw}{F(W)}
\newcommand{\ota}{\ot_A}
\newcommand{\repza}{\rep_\cc^0(A)}
\newcommand{\epsa}{\eps_A}
%%%
\newcommand{\bndefn}{\bn{defn}}
\newcommand{\edefn}{\end{defn}}
\newcommand{\bdefn}{\bn{defn}}

\newcommand{\vld}{V^{*}}
\newcommand{\vldd}{V^{**}}
\newcommand{\xld}{X^{*}}
\newcommand{\xldd}{X^{**}}
\newcommand{\yld}{Y^{*}}
\newcommand{\yldd}{Y^{**}}
\newcommand{\aldu}{A^{*}}
\newcommand{\aldd}{A^{**}}

\newcommand{\ia}{\mtr{i}_A}
\newcommand{\aota}{A\ot A}

\newcommand{\idv}{\id_V}

\newcommand{\ld}{^*}
\newcommand{\repg}{\rep(G)}

\newcommand{\thetav}{\theta_V}

\newcommand{\tta}{\theta_A}

\newcommand{\muv}{\mu_V}
\newcommand{\muw}{\mu_W}

\newcommand{\dimcc}{\dim(\cc)}
\newcommand{\chii}{\chi_i}
\newcommand{\chistar}{\ch_{i^*}}
\newcommand{\chj}{\ch_j}
\newcommand{\chm}{\ch_m}
\newcommand{\chn}{\ch_n}
\newcommand{\dimvi}{\dim(V_i)}%{d(V_i)}%
\newcommand{\mtcd}{Q}%{\mtc D}
\newcommand{\mtca}{\mtc A}
\newcommand{\lamcd}{\lam_\cd}
\newcommand{\fpcd}{\fp(\cd)}
\newcommand{\laml}{\lam_L}
\newcommand{\apm}{A//M}
\newcommand{\apl}{A//L}
\newcommand{\repapm}{\rep(\apm)}
\newcommand{\repapl}{\rep(\apl)}
\newcommand{\dimvj}{\dim(V_j)}%{d(V_j)}
\newcommand{\dvi}{\dim(V_i)}
\newcommand{\dvj}{\dim(V_j)}
\newcommand{\sumjtom}{\sum_{j=0}^m}
\newcommand{\sumitom}{\sum_{i=0}^m}
\newcommand{\sij}{s_{ij}}
\newcommand{\sji}{s_{ji}}
\newcommand{\dxj}{d_j}
\newcommand{\dxi}{d_i}
\newcommand{\dimka}{\dim_{\kk}(A)}
\newcommand{\dimk}{\dim_{\kk}}
\newcommand{\blaml}{\blam_L}
\newcommand{\sumjtor}{\sum_{j=0}^r}
\newcommand{\dimkl}{\dim_{\kk}(L)}
\newcommand{\mtcjl}{\mtc J_L}
\newcommand{\vota}{ V\ot A}
\newcommand{\vi}{V_i}
\newcommand{\vj}{V_j}
\newcommand{\dimcd}{\dim(\cd)}

\newcommand{\alij}{\al_{ij}}
\newcommand{\alji}{\al_{ji}}
\newcommand{\rcc}{r_\cc}
\newcommand{\rcd}{r_\cd}
\newcommand{\clsx}{[X]}
\newcommand{\clsy}{[Y]}
\newcommand{\clsz}{[Z]}
\newcommand{\rcdp}{r_{\cd'}}
\newcommand{\sumjtorp}{\sum_{j=0}^{r'}}
\newcommand{\aljm}{\al_{jm}}
\newcommand{\aljn}{\al_{jn}}
\newcommand{\sjm}{s_{jm}}
\newcommand{\smj}{s_{mj}}
\newcommand{\snj}{s_{nj}}

\newcommand{\betaij}{\beta_{ij}}
\newcommand{\betaji}{\beta_{ji}}

 \newcommand{\ip}{i'}
\newcommand{\sumjtoprp}{\sum_{j=0}^{r'}}
\newcommand{\sumjtopr}{\sum_{j=0}^{r}}
 \newcommand{\teh}{\tilde{h}}
\newcommand{\cdp}{{\cd'}}\newcommand{\xphii}{X_{\phi(i)}}
\newcommand{\inv}{^{-1}}

\newcommand{\fq}{f_Q}
\newcommand{\tr}{\mtr{tr}}
\newcommand{\rtwone}{R_{21}R}

\newcommand{\ccad}{{\cc_{\mtr{ad}}}}
\newcommand{\ccpt}{{\cc_{\mtr{pt}}}}
\newcommand{\qtr}{quasi-triangular\;}
\newcommand{\trq}{\tr_q}

\newcommand{\repal}{\mtr{Rep}(A//L)}
\newcommand{\lkeravi}{\lker_A(V_i)}
\newcommand{\lkeravj}{\lker_A(V_j)}
\newcommand{\cross}[1][1pt]{\ooalign{%
 \rule[1ex]{1ex}{#1}\cr% Horizontal bar
 \hss\rule{#1}{.7em}\hss\cr}}% Vertical bar
\newcommand{\blml}{\blam_L} 
\newcommand{\phir}{\phi_R}
\newcommand{\kda}{{  \Phi(A)}}

\newcommand{\mtcil}{\mtc{I}_L}

\newcommand{\un}{\unu}
\newcommand{\tfl}{\mtc{T}}
\newcommand{\barzm}{\barz(M)}
\newcommand{\barzn}{\barz(N)}
\newcommand{\ccr}{\mtc R^{\cc}}
\newcommand{\ulc}{\ul{\cc}}

\newcommand{\pimx}{\pi_{M;\;X}}
\newcommand{\pinx}{\pi_{N;\;X}}
\newcommand{\acc}{{\mathrm A_\cc}}
\newcommand{\epsu}{\eps_\unu}

\newcommand{\ob}{\mtr{Obj}}
\newcommand{\obc}{\mtr{Obj(\cc)}}
\newcommand{\ccop}{\cc^{\mtr{op}}}
\newcommand{\mtf}{\mtc F_\lam}%{\mtc F}
\newcommand{\mtfi}{\mtc F^{-1}_\lam}
\newcommand{\elcd}{\ell_\cd}
\newcommand{\mcid}{\mtc I_\cd}
\newcommand{\mcidp}{\mtc I_{\cd'}}
\newcommand{\wtildelcd}{\widetilde{\elcd}}
\newcommand{\wtildelcdp}{\widetilde{\ell_{\cd'}}}
\newcommand{\cpt}{\cc_{\mtr{pt}}}
\newcommand{\barzr}{\barz_\cd}
\newcommand{\barzv}{\barz(V)}
\newcommand{\acd}{\mathrm A_\cd}
\newcommand{\czrcd}{\cz_\cc(\cd)}
\newcommand{\sml}{\Small}
\newcommand{\bs}{\blue{\Small }}
\newcommand{\yd}{Yetter-Drinfeld}

\newcommand{\sumitor}{\sum_{i=0}^r}
\newcommand{\cdop}{\cd^{\mtr{op}}}
\newcommand{\ccrev}{\cc^{\mtr{rev}}}
\newcommand{\barz}{{\bar{\mathrm Z}}}
\newcommand{\etl}{etale\;}
\newcommand{\czca}{\cz(\ca)}

%%%%%
\newcommand{\tetx}{\text}
\newcommand{\widehta}{\widehat}
\newcommand{\wdhat}{\widehat}
\newcommand{\wht}{\widehat}
\newcommand{\cofa}{{\mathbb C[A]}}
\newcommand{\wdt}{\widehat}
\newcommand{\dl}{{^\#}}
\newcommand{\comx}{\mathbb C}

\newcommand{\sgj}{\sg(j)}

\newcommand{\mujo}{\mu_\jo}
\newcommand{\mujtw}{\mu_\jtw}
\newcommand{\adz}{a^{\#}}
\newcommand{\bdz}{b^{\#}}

\newcommand{\spr}{S^\perp}
\newcommand{\cofs}{\comp [S]}
\newcommand{\spz}{S^{\perp_z}}

\newcommand{\omz}{\omega_z}
\newcommand{\zg}{\mathrm{Z}(S)}
\newcommand{\aling}{{\al \in g}}

\newcommand{\blkg}{\mtr{Bl}(g)}
\newcommand{\clsg}{\mtr{Cl}(g)}
\newcommand{\mtadinv}{\mtc G^{{-1}}}%^{{\wdht A}}
\newcommand{\muk}{\mu_{k}}
\newcommand{\mta}{\mtc F}%^{A}
\newcommand{\cofad}{\comp[\wdht A]}
\newcommand{\wtau}{\wdht{\tau}}
\newcommand{\mtainv}{{\mta}^{-1}}
\newcommand{\wdht}{\widehat}
\newcommand{\augm}{\mtr{aug}}
\newcommand{\mua}{\wdht {\wdht a}}
\newcommand{\aps}{A//S}
\newcommand{\cssa}{\cc(S, A)}
\newcommand{\aug}{\mtr{aug}}
\newcommand{\rss}{{\big|_S}}
\newcommand{\gprp}{g^\perp}
\newcommand{\alins}{{s \in S}}

\newcommand{\sz}{s^{D}}
\newcommand{\wmu}{\widehta{\mu}}
\newcommand{\wmui}{\widehta{\mu}_i}
\newcommand{\wmuj}{\widehta{\mu}_j}
\newcommand{\wch}{\widehta{\ch}}

\newcommand{\wzd}{\widehat{d}}
\newcommand{\wpm}{\widehat{P}}
\newcommand{\wps}{\widehat{p}}

\newcommand{\gal}{\mtr{Gal}}
\newcommand{\galkq}{\gal(\mathbb K/\mathbb Q)}
\newcommand{\sgf}{\sg_{_F}}
\newcommand{\sggi}{{\sg(i)}}
\newcommand{\sge}{\sg_{_E}}
\newcommand{\unue}{{\unu_{\cecc}}}

\newcommand{\mtcf}{\mtc {F}}

\newcommand{\wsgf}{\widehat{{\sg}_{F}}}

\newcommand{\we}{\widehta{E}}
\newcommand{\sumktom}{\sum_{k=0}^m}

\newcommand{\wf}{\widehat{F}}

\newcommand{\hsgj}{\widehat{\sg}(j)}
\newcommand{\whsgi}{\widehta{\sg}(i)}

\newcommand{\wpp}{\widehat{p}}
\newcommand{\tauj}{{\tau(j)}}
\newcommand{\dimcctauj}{\dim(\cc^\tauj)}
\newcommand{\etas}{{\eta(s)}}
\newcommand{\mcc}{m_\cc}

\newcommand{\wal}{\widehta{\al}}
\newcommand{\wj}{\widehat{\mtc J}}
\newcommand{\galc}{\mtr{Gal}_{\cc}}
\newcommand{\galz}{\mtr{Gal}_{\czcc}}
\newcommand{\wjr}{\widehat{J}_{R}}

\newcommand{\dimcck}{\dim(\cc^k)}

\newcommand{\wgrcc}{\widehat{\mtr{Gr}(\cc)}}
\newcommand{\nchi}{{\frac{\ch_i}{d_i}}} \newcommand{\nchj}{{\frac{\ch_j}{d_j}}}
\newcommand{\wni}{\widehat{n_i}}

\newcommand{\sgte}{\widetilde{\sg_E}}

\newcommand{\mtad}{\mtc G}
\newcommand{\whj}{\widehta{h}_j}
\newcommand{\jdl}{{j\dl}}
\newcommand{\wcfcc}{\widehat{\cfcc}}
\newcommand{\mutauj}{\mu_{\tau(j)}}
\newcommand{\tauk}{\tau(k)}
\newcommand{\muzm}{{\mu_0^{-}}}
\newcommand{\sqrtog}{\sqrt{|G|}}
\newcommand{\muz}{\mu_0}
\newcommand{\njtw}{n_\jtw}
\newcommand{\njo}{n_\jo}
\newcommand{\fjo}{F_\jo}
\newcommand{\fjtw}{F_\jtw}
\newcommand{\wta}{\widehat{A}}

\newcommand{\dol}{{^{\circ}}}
\newcommand{\bdl}{{b\dl}}
\newcommand{\jdol}{{j\dol}}
\newcommand{\fj}{F_j}

\newcommand{\cwta}{\comp[\wta]}

\newcommand{\hx}{\widehta{x}}
\newcommand{\hy}{\widehta{y}}

\newcommand{\cal}{\mtc A_{\al}}
\newcommand{\xuu}{x_{uu}}
\newcommand{\wxuu}{\widehat{\xuu}}
\newcommand{\xvv}{x_{vv}}
\newcommand{\xuv}{x_{uv}}
\newcommand{\xmn}{x_{m,n}}
\newcommand{\buvmn}{B^{u,v}_{m,n}}
\newcommand{\blm}{\blam}
\newcommand{\dimccr}{\dim(\cc^r)}
\newcommand{\adl}{a\dl}
\newcommand{\sumltom}{\sum_{l=0}^m}

\newcommand{\mbq}{\mathbb Q}
\newcommand{\mbqs}{\mathbb Q(S)}
\newcommand{\mbk}{\mathbb K}
\newcommand{\mz}{\mathbb Z}

\newcommand{\wsgj}{\widehat{\sigma}(j)}
\newcommand{\wsgi}{\widehat{\sigma}(i)}
\newcommand{\wg}{\widehat{g}}
\newcommand{\wtf}{\widehat{F}}
\newcommand{\galqspq}{\mtr{Gal}(\mathbb Q(S)/\mathbb Q)}
\newcommand{\cctauj}{\cc^{\tau(j)}}
\newcommand{\cctauk}{\cc^{\tau(k)}}
\bd
\newcommand{\wtfj}{\widetilde{F_j}}
\newcommand{\wfj}{\widetilde{F_j}}
\newcommand{\wtmuj}{\widetilde{\mu_j}}
\newcommand{\wmtcfj}{{\widetilde{\mtc F}_j}}
\newcommand{\mtfr}{\mtr{F_a}}
\newcommand{\wdr}{\widehat {R_\comp}}

\newcommand{\fgph}{{F_{G/H}}}
\newcommand{\wcfj}{\wmtcfj}

\newcommand{\nxi}{{\frac{x_i}{d_i}}}
\newcommand{\fpr}{{\fp(R)}}
\newcommand{\nxs}{{\frac{x_s}{d_s}}}

 \newcommand{\mtfme}{\mtc F}
\newcommand{\chic}{\ch_i^{\circ}}
\newcommand{\chjc}{\ch_j^{\circ}}
\newcommand{\mtfsh}{{\mtc F_\lam}}%{\ul{\mtc F}}
\newcommand{\mupq}{{\mu_{pq}}}

\newcommand{\tlam}{{\widetilde{\lam}}}
\newcommand{\chid}{{\ch_i^{\circ}}}
\newcommand{\rc}{{R_\comp}}
\newcommand{\rgo}{{\mathbb R_{\geq 0}}}

\newcommand{\aliro}{{\al_{i\ro}}}
\newcommand{\alm}{\lam}
\newcommand{\mtrs}{{\mtr S}}
\newcommand{\sumroirc}{{\sum_{\ro \in \irr(\rc)}}}
\blue{\title[Isaacs property and Extended Haagerup]
{Structure constants, Isaacs property and Extended Haagerup fusion categories}}

\author{Sebastian Burciu}
\address{Inst.\ of Math.\ ``Simion Stoilow" of the Romanian Academy P.O. Box 1-764, RO-014700, Bucharest, Romania}
\email{sebastian.burciu@imar.ro}

\author{Sebastien Palcoux}
\address{S. Palcoux, Yanqi Lake Beijing Institute of Mathematical Sciences and Applications, Huairou District, Beijing, China}
\email{sebastien.palcoux@gmail.com}
\urladdr{https://sites.google.com/view/sebastienpalcoux}

\thanks{The first author is supported by a grant of the Ministry of Research, Innovation and Digitization, CNCS/CCCDI - UEFISCDI, project number PN-III-P4-ID-PCE-2020-0878, within PNCDI III. The second author is supported by BIMSA Start-up Research Fund and Foreign Youth Talent Program from the Ministry of Sciences and Technology of China (Grant QN2021001001L)}
%\date{\today}

\maketitle

\begin{abstract}
This paper presents an abstract Isaacs property that involves the Fourier transform for fusion rings, which may be non-commutative, thus expanding upon the commutative version described in \cite{lpr1}. A categorical version of this property was subsequently introduced in \cite{eno-nec} for any spherical fusion category, matching with our abstract version in the pseudo-unitary case. We demonstrate that the Isaacs property occupies a distinct position, falling between the integrality of structure constants and the $1$-Frobenius properties, in the commutative case. We show that the Extended Haagerup fusion categories, denoted as $\mathcal{EH}_i$, do not satisfy the Isaacs property. This finding provides a negative response to \cite[Question 5.8]{eno-nec}, refutes \cite[Conjecture 2.5]{lpr1}, and recovers that $\mathcal{EH}_1$ lacks a braiding structure.
\end{abstract}

\section{Introduction}
Recently, various criteria for a fusion ring to be categorifiable were intensively studied, see \cite{eno-nec, lpw, lpr1, lpr2, lpr3, hlpw} and the references therein. Following \cite{lpr1}, a categorical Isaacs property was formulated in \cite{eno-nec} for any fusion category, and \cite[Proposition 5.2]{eno-nec} states that every braided spherical fusion category satisfies this property.

We introduce an abstract Isaacs property involving the Fourier transform for (possibly non-commutative) fusion rings, extending the (commutative) one introduced in \cite{lpr1}, and matching with the one in \cite{eno-nec} in the pseudo-unitary case.

We remark that the version in \cite{eno-nec} involves the dimension of some simple objects of the Drinfeld center $\czcc$. However, we show that this property can be formulated purely algebraically in terms of the ring structure of the Grothendieck ring $K(\cc)$, if $\cc$ is pseudo-unitary.

It is conjectured in \cite{eno-weakly} that every fusion category $\cc$ is $1$-Frobenius, i.e the ratio  $\frac{\fp(\cc)}{\fp(X)}\in \mathbb A$ (the ring of algebraic integers) for every simple object $X$ of $\cc$. Moreover, \cite[Proposition 5.4]{eno-nec} states that a fusion category with the categorical Isaacs property and a commutative Grothendieck ring is $1$-Frobenius. Theorem \ref{is-abstr} extends this result for arbitrary fusion rings with Isaacs property but with the additional hypothesis that the basic element is central.

Structure constants for pivotal fusion categories with a commutative Grothendieck ring were introduced in \cite{b-blms}. In this paper we extend this notion to non-commutative fusion rings using the {\it matrix class sums} coming from the central primitive idempotents of the fusion ring.  Using these matrix-class sums, we prove a Frobenius divisibility type result for commutative fusion rings in Theorem \ref{d-Frobenius}.

We finally prove that the Extended Haagerup fusion categories, introduced in \cite{gmpps} and denoted $(\mathcal{EH}_i)_{i=1,\dots,4}$, are not Isaacs, recovering that $\mathcal{EH}_1$ (the only one with a commutative Grothendieck ring) has no braiding (first proved in \cite{MoWa}), providing a negative answer to \cite[Question 5.8]{eno-nec} asking whether every spherical fusion category satisfies the categorical Isaacs property, and disproving \cite[Conjecture 2.5]{lpr1}. Moreover the Extended Haagerup fusion categories are the only simple fusion categories known to be non-Isaacs, but we can make infinitely many (non-simple) ones by Deligne tensor product.

We provide a sufficient condition (involving the Morita equivalence) for a property to be true for every spherical fusion category (see \S \ref{sub:MoEq}). We deduce that the $1$-Frobenius property holds for every spherical fusion category if and only if it is invariant by Morita equivalence. Idem for the Isaacs property, so that it cannot be invariant by Morita equivalence, as the  Extended Haagerup fusion categories are not Isaacs.

Theorem \ref{d-Frobenius} implies that the Isaacs property is positioned between the integrality of the structure constants and the $1$-Frobenius properties (in the commutative case), and it is strictly in between thanks to $\mathcal{Z}(\mathrm{Vec}_{S_3})$ on one hand (see \cite{CW6}), and $\mathcal{EH}_1$ on the other hand. The fact that the Isaacs property fails on exotic known examples only suggests the possibility of even more exotic counterexamples to the $1$-Frobenius property.

Shortly, this paper is organized as follows: \S \ref{afr} studies  basic properties of fusion rings; \S \ref{gfc} introduces the matrix class sums and study their basic properties; \S \ref{is} introduces the abstract Isaacs property for fusion rings and prove that it coincides with the categorical Isaacs property from \cite{eno-nec}. In this section we also prove the two Frobenius divisibility type results mentioned above. Finally, \S \ref{sec:EH} is dedicated to show that the Extended Haagerup fusion categories are not Isaacs.

\tableofcontents

\section{Preliminaries on abstract fusion rings}\label{afr}
Recall \cite[\S 3]{EGNO15} that a {\it fusion ring} $(R, \mtc B)$ is a ring $R$ which is free as $\mathbb Z$-module  with a finite basis $\mtc B=\{x_0,x_1,\dots,x_m\}$, called {\it standard basis}, satisfying the following properties:
\bne
\item $x_0=1$ is the unit of $R$,
\item $x_ix_j=\sumktom N_{ij}^k x_k$ with $N^k_{ij}\in \mathbb Z_+$,
\item there is  ring involution $^*$ on $R$ which is stable on $\mtc B$ and such that $N^0_{ij}=\delta_{i,j^*}$ (where $i^*$ is defined by $x_{i^*}:=x_i^*$). 
\ene 
%for all $i$ there is $i^*$ such that $N^0_{ij}=\delta_{i,j^*}$, inducing an involution on $\mtc B$ defined by $x_i^*:=x_{i^*}$.
\subsection{The trace $\tau$ and its non-degenerate associative bilinear form.} \label{taur}  %
The involution on the basis $\mtc B$ induces a $*$-structure on the finite dimensional algebra $R_\comp:=R\ot_{\mathbb Z}\comp$ making it a semisimple algebra. By Wedderburn-Artin theorem one has:
$$
\rc\simeq \prod_{\rho\in \irr(\rc)} M_{\deg \rho}(\comp)
$$
Since $\rc$ is a semisimple $\comp$-algebra, by abuse of notations, in this paper we identify the irreducible representations of $\rc$  with their characters. Recall that one can define a linear function $\tau:\rc\ra \comp$ with $\tau(x_i)=\delta_{i,0}$, where as above, $x_0=1$. By results covered in \cite[\S 3]{EGNO15}, it follows that  $\tau:\rc\ra \comp$ is a trace, i.e $\tau(ab)=\tau(ba)$, for all $a, b\in \rc$. Moreover the bilinear form $(\;,\;)_{\tau}:\rc\otimes \rc\ra \comp, (a,b)\sent \tau(ab)$ is associative symmetric non-degenerate and therefore one can write
\beq\label{trform}
(\;,\;)_{\tau}:=\sum_{\rho\in \irr(\rc)}\frac{1}{n_\rho}Tr_\rho.
\eeq
for some non-zero scalars $n_\rho \in \comp^\times$.
Since  $\{x_i,x_{i^*}\}$ is a pair of dual bases for $(\;,\;)_{\tau}$ it follows that %\blue{usual reference?? curtis and reiner-maybe}
\beq\label{db}
\sum_{\rho \in \irr(\rc)} \sum_{p,q=1}^{\deg \ro}F^\rho_{pq}\ot {n_\rho}F^\rho_{qp}=\sumitom x_i\ot x_{i^*}.
\eeq
where $\{F^\rho_{pq}\}_{1\leq p,q \leq \deg \rho}$ is a linear matrix-basis for the block $M_{\deg \rho}(\comp)$. Note that one has
$$
f_\rho:=\sum_{x\in \mtc B}\rho(x)x^*=\sum_{1\leq p,q\leq \deg(\rho)}n_\rho\rho(F^\rho_{pq})F^\rho_{qp}=n_\rho F^\rho 
$$
where $F^\rho:=\sum_p F^\rho_{pp}$ is the central primitive idempotent of $\rc$ corresponding to $\rho \in \irr(\rc)$. Therefore, as in \cite{o-fc}, one can define the {\it formal codegree} $c_\rho$ of $R$ at $\rho$,  as  the scalar by which $f_\rho$ acts on $\rho$. Thus, with the above notations, one has that the formal codegree  $c_\rho$ equals the scalar $n_\rho.$
\subsection{A multiplication on $\wdr$} 
Let $(R, \mtc B)$ be a fusion ring.   For any element $x_i \in \mtc B$ denote $d_i:=\fp(x_i)$, the Frobenius-Perron dimension of $x_i$. Denote also by $\wdr$ the linear dual of $\rc$, i.e. $\wdr:=(\rc)^*.$
Following \cite{b-blms, hdk} one can define a  multiplication on $\wdr$ in the following way. For any $\mu, \nu \in \wdr$, the linear map  $\mu\star \nu\in\wdr $ is defined on the basis $\{\nxs\}$ by: 
\beq\label{mulpsi}
[\mu\star \nu](\nxs):=\mu(\nxs)\nu(\nxs).
\eeq
Then $\mu\star \nu$ is linearly extended on the whole $\rc$. Clearly $\wdr$ becomes a  commutative algebra. 
\begin{notations}  
We denote by $\{\rho_{pq}\in \wdr\}_{\rho\in \irr(\rc),\;0\leq p,q\leq \deg\rho}$ the linear dual basis of the matrix-basis $\{F^\rho_{pq}\}$ of $\rc$. Therefore  for any two irreducible representations $\rho, \psi\in \irr(\rc)$ one has $\rho_{pq}(F^\psi_{p'q'})=\delta_{\psi, \rho}\delta_{p,p'}\delta_{q,q'}$. Denote also by $\{x_i^\circ\in \wdr\}_{i=0}^m$ the linear dual basis of $\{x_i\}$. Therefore $\lag x_i^\circ, x_j\rag=\delta_{i,j}$ and note that $x_0^\circ=\tau.$
\end{notations}

\bl\label{didempt} 
For any $0\leq i \leq m$ one has that $\widetilde E_i:=d_ix_i^\circ$ are the orthogonal primitive idempotents of $\wdr$. The linear character $\tomega_i:\wdr \ra \comp$ corresponding to $\widetilde E_i$ is given by 
$$
\tomega_i:\wdr\ra \comp,\;\;\mu\sent \frac{\mu(x_i)}{d_i}.
$$
\el
\bpf
Indeed, note that
$$[x_i^\circ\star x_j^\circ](\nxs)=x_i^\circ(\nxs)x_j^\circ(\nxs)=\delta_{i,s}\delta_{j,s}\frac{1}{d_s^2}=\delta_{i,j}\frac{1}{d_i}x_i^\circ(\nxs).
$$
On the other hand one has
$
\frac{1}{d_i}x_i^{\circ}(\nxs)=\delta_{i,s}\frac{1}{d_s^2}
$, i.e  $x_i^\circ\star x_j^\circ=\delta_{i,j}\frac{1}{d_i}x_i^{\circ}$.

Note also
\begin{eqnarray*}
[\widetilde E_i\star \nu](\frac{x_s}{d_s})&=&[d_ix_i^\circ\star\nu](\frac{x_s}{d_s})=d_i{x_i^\circ}(\frac{x_s}{d_s})\nu(\frac{x_s}{d_s})\\&=&\delta_{s,i}\nu(\frac{x_s}{d_s})=\widetilde E_i(\frac{x_s}{d_s})\nu(\frac{x_i}{d_i})= \widetilde E_i(\frac{x_s}{d_s})\tomega_i(\nu)
\end{eqnarray*}
which shows that
$
\widetilde E_i\star \nu=\tomega_i(\nu)\widetilde E_i
$
i.e.  $\tomega_i(\mu):=\mu(\frac{x_i}{d_i})$ are the characters of $\wdr$.
\epf

\subsection{A Fourier transform}\label{fouriertr}
Define a $\comp$-linear map $\mtc F: \rc\ra \wdr, \;x_i\sent \fpr x_{i^*}^\circ.$
Clearly $\mtc F$ is bijective and on the linear basis $\{x_{i^*}^\circ\}$ its inverse is given by
\beq\label{finvlin}
\mtc F^{-1}(x_i^\circ)=\frac{1}{\fpr}x_{i^*}.
\eeq

Recall that $\fp:R\ra \comp$ is a linear character of $\rc$. Denote by $F_0:=F^{\fp}\in R$ the primitive central idempotent associated to $\rho=\fp$. Next we show that the inverse of the Fourier transform $\mtc F$ is related to the following functional
$$
\mtc G:\wdr\ra \rc, \mu\sent (\nu\sent \lag\mu\star\nu, F_0\rag).
$$ 
Note that above $\mtc G$ is defined by using the usual duality $\widehat{\wdr}\simeq \rc$.
\bp 
Let $(R, \mtc B)$ be a fusion ring.  With the above notations one has \beq\label{mtfinv1}
\mtc G \circ \mtc F=(-)^*.
\eeq
\ep
\bpf
It is enough to show that
$\mtc G\circ \mtc F(x_i)=x_{i^*}$, i.e 
$$
\mtc G(\fpr x_{i^*}^{\circ})=x_{i^*}.
$$
Applying $\fp \ot \fp$ to Equation \eqref{db}, we get $n_{\fp} = \fp(R)$, next applying $\fp \ot \id$, we get
\beq\label{fzero}
F_0=\frac{1}{\fpr}\big(\sum_{x_i\in B}d_{i}x_i\big).
\eeq
Note that $d_i = \fp(x_i)=\fp(x_i^*) = d_{i^*}$.  By the definition of $\mtc G(\mu)$, for any $\mu\in \wdr$ one has that
\begin{eqnarray*}
\lag\nu, \mtc G(\mu)\rag&=&\lag\mu\star \nu, F_0\rag\\&=&\lag\mu\star \nu, \frac{1}{\fpr}\bigg(\sum_{x_j\in B}d_{j} x_j\bigg)\rag 
\\&=& \frac{1}{\fpr}\bigg(\sum_{x_j\in B} d_j^2\lag\mu\star \nu, \;\frac{x_j}{d_j}\rag\bigg)
\\&=& \frac{1}{\fpr}\bigg(\sum_{x_j\in B} d_j^2\lag\mu, \;\frac{x_j}{d_j}\rag\lag\nu, \;\frac{x_j}{d_j}\rag\bigg)
\end{eqnarray*}
Therefore 
\beq\label{nugmu}
\lag\nu, \mtc G(\mu)\rag=\frac{1}{\fpr}\bigg(\sum_{x_j\in B} \lag\mu, \;{x_j}\rag\lag\nu, \;{x_j}\rag\bigg)
\eeq
For $\mu=\mtc F(x_i)=\fpr x_{i^*}^\circ$ one has 
\begin{eqnarray*}
 \lag\nu, \mtc G(\fpr x_{i^*}^\circ)\rag&=&\fpr \frac{1}{\fpr}\lag\nu, \;{x_{i^*}}\rag=\\&=&\lag\nu, \; {x_{i^*}}\rag
\end{eqnarray*}
Since the last equality holds for any $\nu\in \wdr$, we get $\mtc G(\mtc F(x_i))=\mtc G(\fpr x_{i^*}^\circ)={x_{i^*}}$, and so Equation \eqref{mtfinv1} holds.
\epf
Equation \eqref{mtfinv1} implies that 
$\mtc G \circ \mtc F\circ (-)^*= \id_{R_\comp}, \;\text{and}\; (-)^*\circ \mtc G \circ \mtc F= \id_{R_\comp},$
which gives that
\beq\label{gfinv}
\mtc G^{-1}=\mtc F\circ (-)^*, \;\text{and}\;\; \mtc F^{-1}=(-)^*\circ \mtc G.
\eeq
\section{On Grothendieck rings of pivotal  fusion categories}\label{gfc}
Let $\cc$ be a fusion category and $R:=K_0(\cc)$ its Grothendieck ring. Let $\irr(\cc):=\{X_0, X_1, \cdots , X_m\}$ be a complete set of isomorphism representatives for the simple objects of $\cc$. It is well known that $R$ is a fusion ring with standard basis $\{[X_i]\}_{i=0}^m$, where $[X]$ stand for the isomorphism class of the object $X$. Therefore all the results of the previous section can be applied in these settings for $R=K_0(\cc)$. Let $K(\cc):=R_\comp=K_0(\cc)\ot_{\mathbb Z}\comp$ be the complex Grothendieck ring of $\cc$.

Let $\czcc$ be the Drinfeld double of $\cc$ and $F:\czcc\ra \cc$ the forgetful functor. Then $F$ admits a right adjoint $R$ and $Z:=FR:\cc \ra \cc$ is a Hopf comonad , called the {\it central Hopf comonad} associated to $\cc$, see \cite[\S 3.1]{scalg}. 

It is well known that $A:=Z(\unu)$ has the structure of central commutative algebra in $\cz(\cc)$.

The vector space $\cecc:= \hm_{\C}(\unu, A) $ is called {\it the space of central elements.} On this space one can define a multiplication  such that  $z.w=m\circ (z\ot w)$ where $m:A\ot A\ra A$ is the multiplication of the central commutative algebra $A$.

The vector space $\cfcc:=\hm_\cc(A, \unu)$ is called the {\it space of class functions} of $\cc$.  For two class functions $f, g\in \cfcc$ one can also  define a multiplication by
 $f\star g:=f \circ Z(g) \circ \delta_{\unu},$
where $\delta: Z \ra  Z^2$ is the comultiplication structure of the Hopf comonad $Z$,  see \cite{scalg}. 

Let now $\cc$ be a pivotal fusion category with the pivotal structure denoted by  $j:\id_\cc\ra (-)^{**}$. For any object $X$ of $\cc$, with the help of the pivotal structure $j$ Shimizu has defined in \cite{scalg} a class function $\mtr{ch}(X)\in \cfcc$. 

By \cite[Theorem 3.10]{scalg} one has that $\mtr{ch}(X\ot Y)=\mtr{ch}(X)\mtr{ch}(Y)$ for any two objects $X$ and $Y$ of $\cc$. It was also shown in \cite[\S 4]{scalg} that the function  $K(\cc)\ra \cfcc,\;[X]\ra \mtr{ch}(X)$ is an isomorphism of $\comp$-algebras. Since $K(\cc)$ is a semisimple algebra one can write a Wedderburn-Artin decomposition by
\beq\label{wdb}
\cfcc\simeq \bigoplus_{\ro\in \irr(\cfcc)}M_{\deg\rho}(\comp).
\eeq

Shimizu in \cite{scalg},  defined a non-degenerate pairing
$$\langle\;,\; \rangle_z : \cfcc \times \cecc\ra \unu,$$ given by $ \langle f, a\rangle_z \id_{\unu}= f \circ a,$
for all $f \in \cfcc$ and $a\in \cecc$. 

Recall $R:\cc \ra \czcc$ is a right adjoint to the forgetful functor $F:\czcc \ra \cc$. As explained in \cite[Theorem 3.8]{scalg} this adjunction  gives an isomorphism algebras
\beq\label{adjisom}
\cfcc \xra{\cong} \mtr{End}_{\czcc}(R(\unu)),\;\; \ch\mapsto Z (\ch)\circ \delta_\unu.
\eeq

The above isomorphism combined with Equation \eqref{wdb} allows us  to write $R(\unu)=\bigoplus_{\ro\in \irr(\cfcc)}\ccro$ for the decomposition of $R(\unu)$ in homogeneous components in $\czcc$. Note that each homogeneous component can be written as
$\ccro=\bigoplus_{s=1}^{\deg \rho}\ccpsro$
where $\ccpsro$ are the simple (isomorphic) sub-objects of $R(\unu)$ entering in the homogeneous component $\ccro$. Therefore as an object of $\czcc$ one has a decomposition in simple objects
\beq\label{ss}
R(\unu)=\bigoplus_{\ro\in \irr(\rc)}\bigoplus_{1\leq s\leq \deg\rho}\ccpsro.
\eeq
 Following \cite{scalg} a  {\it cointegral} in $\cc$  is the unique element (up to a scalar) $  \lambda\in \cfcc $  such that $\ch_i\lam=\lam\ch_i=\dim(X_i)\lam$ for any irreducible character  $\ch_i:=\mtr{ch}(X_i)$. Here $\dim(X_i)$ is the categorical dimension of $X_i$.

Furthermore, let as above, $\irr(\cc):=\{X_0, \dots, X_m\}$ be a complete set of representatives of isomorphism classes of simple objects. As in previous section, let $d_i:=\fp(X_i)$ the Frobenius-Perron dimension of $X_i$. To any simple object $X_i$ of $\cc$ Shimizu has associated in  \cite{scalg} the corresponding primitive central elements $E_i\in \cecc$ such that $\lag\ch_i, E_j\rag_z=\dim(X_i)\delta_{i, j}$ where $\ch_i:=\mtr{ch}(X_i)$ is the irreducible character associated to the simple object $X_i$. One has that $\{E_i\}_{i=0}^m$ form a linear basis of $\cecc$ and $E_i.E_j=\delta_{i,j}$.

Without loss of generality we may suppose that $X_0=\unu$. It is easy to see that in this case $\ch_0=\epsu$ is the the counit of the Hopf comonad $Z$ and unit of the algebra $\cfcc$.  

For any $i \in \{0, \dots, m\}$, we define $i^*\in \{0, \dots, m\}$ by $X_i^*\simeq X_{i^*}$. Then $i \sent i^*$ is an involution on $\{0, \dots, m\}$.
By \cite[Equation 6.8]{scalg} one has that the idempotent cointegral of $\cc$ has the form
\beq\label{intregform}\lam_{\cc}=\frac{1}{\dimcc}(\sum_{[X_i]\in \irr(\cc)}\dim(X_i^*)\ch_{i}).
\eeq 
\subsection{Dual $\wcfcc$ of the Grothendieck ring}
For  $R=K_0(\cc)$ denote the corresponding trace $\tau_\cc:=\tau$.  Then the symmetric associative non-degenerate bilinear form on $K(\cc)\simeq \cfcc$ is given by $(\ch, \mu)_\cc:=\tau_\cc(\ch\mu)$. Suppose as above that $
\tau_\cc=\sum_{\psi\in\irr(\cfcc)}  \frac{1}{n_\psi}Tr_\psi$ for some non-zero scalars $n_\psi$. 

As in the previous section, since $\cfcc$ is a semisimple algebra, one can  write
$$
\cfcc\simeq \prod_{\rho\in \irr(\cfcc)} M_{\deg \rho}(\comp).
$$
Recall that, as in the previous section  we may fix a linear matrix-basis $\{F^\ro_{pq}\}$ of $\cfcc$ consisting  of the entries of each matrix block $M_{\deg \rho}(\comp)$.  As previously, also denote  by $\{\rho_{pq}\in \wcfcc\}_{\rho\in \irr(\cfcc),\;0\leq p,q\leq \deg\rho}$ the linear dual basis of this matrix-basis. Therefore $\rho_{pq}(F^\psi_{p'q'})=\delta_{\psi, \rho}\delta_{p,p'}\delta_{q,q'}$ for any  $\rho, \psi\in \irr(\cfcc)$.

By  \cite[Lemma 3.27]{b-etale} it follows that in a pivotal fusion category  one has
\beq\label{nro-etale}
{n_\rho}=\frac{\dim(\cc)}{\dim(\cc_1^\rho)}=\frac{\deg\rho\dim(\cc)}{\dim(\cc^\rho)}
\eeq
where $\cc_1^\rho$ is a simple object of the homogeneous component $\cc^\rho$.

Define $\wcfcc$ as the linear dual vector space of $\cfcc$. Clearly, as in the previous section, this is a commutative algebra with multiplication:
$$[\mu\star \nu](\frac{\ch_i}{d_i})=\mu(\frac{\ch_i}{d_i})\nu(\frac{\ch_i}{d_i}),
$$
for all $\mu, \nu \in \wcfcc$. By Lemma \ref{didempt} one has as above that $\widetilde E_i:=d_i\ch_i^\circ\in \wcfcc$ are the orthogonal primitive idempotents of $\wcfcc$ and $\widetilde{\omega}_i(\mu):=\mu(\nchi)$ are the corresponding irreducible characters of the dual $\wcfcc$. 

Let  $\lambda\in \cfcc$ be  the non-zero idempotent cointegral of $\cc$. Shimizu introduced a {\it Fourier transform} of $\cc$  associated to $\lambda$ as the linear map
\beq
\mtfsh:\cecc\ra \cfcc\;\;\text{given by}\;\;a \mapsto \lambda \lh \mtc S(a).
\eeq
 Shimizu has also shown in \cite[Equation (6.10)]{scalg} that 
\beq\label{shm6.10}
\mtfsh(E_i)=\frac{\dim(X_i)}{\dim(\cc)}\ch_{i^*},\; \mtfsh^{-1}(\ch_i)=\frac{\dim(\cc)}{\dim(X_{i^*})}E_{i^*},
\eeq
Note that in \cite{scalg} the author  assumed that the Grothendieck ring of $\cc$ is commutative but his proof for Equation \eqref{shm6.10} works also in the general case. Therefore by \cite[Equation (4.7)]{ccc-march} one has:
\beq\label{myeq}
\lag\ch,\;\mtfsh^{-1}(\mu) \rag_z=\dim(\cc)\tau(\ch\mu),
\eeq
 for all $\ch ,\mu \in \cfcc$ (see also \cite[Equation (2.17)]{b-etale}.)
Then from the definition of $\tau$  this implies 
\begin{eqnarray*}
\lag \ch, \;\mtfsh^{-1}(F^{\ro}_{pq})\rag_z&=&\dim(\cc)\tau(\ch F^{\ro}_{pq})=\dimcc \ro_{qp}(\ch)\tau(F^\ro_{qq})\\&=&\ro_{qp}(\ch)\frac{\dim(\cc)}{n_\ro}
\end{eqnarray*}
for any $\ch\in \cfcc$. Equation \eqref{nro-etale} gives now
\beq\label{eval}
\frac{1}{\dim(\cc_1^\ro)}\lag\ch_i,\;\mtfsh^{-1}(F^\ro_{pq})\rag_z=\ro_{qp}(\ch_i).
\eeq
Define the {\it matrix class sum} of $\cc$ associated to $\rho \in \irr(R_\comp)$ as the central element $\mtr C^{\rho}:=\mtfsh^{-1}(F^{\rho})\in \cecc.$ For pivotal fusion categories with a commutative Grothendieck ring this notion was previously introduced in \cite{b-blms}. Note that by Equation \eqref{eval}, since $\ro=\sum_{p}\ro_{pp}$ one has that
\begin{eqnarray*}
\ro(\ch_i)&=&\sum_p\ro_{pp}(\ch_i)=\sum_p\frac{\lag\ch_i,\;\mtc F_\lam^{-1}(F^\ro_{pp})\rag_z}{\dim(\cc_1^\ro)}=\\&=&\frac{\lag\ch_i,\;\mtc \sum_pF_\lam^{-1}(F^\ro_{pp})\rag_z}{\dim(\cc_1^\ro)}=\frac{\lag\ch_i,\;\mtr{C}^\ro\rag_z}{\dim(\cc_1^\ro)}.
\end{eqnarray*}
Thus
\beq\label{mujev}
\ro(\ch_i)=\frac{\lag\ch_i,\;\mtr{C}^\ro\rag_z}{\dim(\cc_1^\ro)}.
\eeq 
\subsection{On the canonical isomorphism $\al$} Let $\cc$ be a pivotal fusion category. With the above notations remark that both $\wcfcc$ and $\cecc$ are commutative $\comp$-algebras of dimension equal to the rank of $\cc$. In this subsection  a canonical isomorphism $\al:\wcfcc\ra\cecc, \widetilde E_i\sent E_i$ between these two algebras is constructed. In the case of a pivotal fusion category with a commutative Grothendieck ring this isomorphism $\alpha$ was constructed in \cite[Theorem 3.4]{b-blms}. Note the much simpler description of $\al$ given here in terms of the primitive central idempotents of both algebras. Then Equation \eqref{eval} shows that $\al$ coincides to the isomorphism constructed in \cite[Theorem 3.4]{b-blms} in the case of a commutative Grothendieck ring. 

It is also not difficult to check that  $\al$ is the unique linear isomorphism $\beta:\wcfcc\ra\cecc$ such that $\lag \ch,\;\beta(\ro)\rag_z=\ro(\ch)$, for any $\ch \in \cfcc$ and $\ro\in \wcfcc$

Note that \cite[Lemma 30]{b-jpaa} shows that  the linear map $\omega_i$ defined by $ z\sent \frac{\ch_i(z)}{\dim(X_i)}$ is a linear character of the space of central elements. We  call $\omega_i$ the {\it central character} associated to $X_i$. 
By the definition of $\al$ then clearly
\beq\label{moveomega}
\omega_i \circ \al=\widetilde{\omega_i},
\eeq
where $\widetilde{\omega_i}$ is defined in Lemma \ref{didempt}. As in \S \ref{fouriertr}, associated to the fusion ring $R=K_0(\cc)$, one can also define a Fourier transform $\mtfme: \cfcc\ra \widehat{\cfcc} $ given by $\ch_i\sent  \fp(\cc)\ch_{i^*}^{\circ}.$

\bp
With the above notations, if $\cc$ is a pseudo-unitary category one has that 
\beq\label{msh}
\mtfsh^{-1}=\al \circ \mtc F.
\eeq
\ep
\bpf

By the definition of $\al$ one has $ \al(\ch_i^\circ)=\al(\frac{\widetilde E_i}{\fp(X_i)})=\frac{E_i}{\fp(X_i)}$. 
It follows that $\al(\mtc F(\ch_i))=\dim(\cc)\al(\ch_{i^*}^\circ)=  \frac{\dimcc}{\fp(X_{i^*})}E_{i^*}$. On the other hand, by Equation \eqref{shm6.10} one has 
$ \mtfsh^{-1}(\ch_i)=\frac{\dim(\cc)}{\dim(X_{i^*})}E_{i^*}$ which proves the desired equality $\al \circ \mtc F=\mtfsh^{-1}$ since in the pseudo-unitary case $\fp(X_i)=\dim(X_i)$.
\epf
\section{Isaacs property}\label{is}
In this section we show that the Isaacs property for a pseudo-unitary fusion category $\cc$ as defined in \cite{eno-nec} is actually a property of the  Grothendieck ring $K_0(\cc)$ not of the category $\cc$.
\subsection{Isaacs property for fusion rings}
Let $(R, \mtc B)$ be a fusion ring. With the notations from \S \ref{afr}, one can define abstract {\it matrix class sums} by $\mtr S^\ro:={\mtc F}(F^\ro)\in \wdr,$ for any $\rho\in \irr(\rc)$.
\bp\label{mstep}
With the above notations one has
\beq\label{3termeq}
\tomega_i(\mtr S^\ro)=\frac{\fpr}{\fp(x_i) } \frac{\ro(x_i)}{c_\ro}.
\eeq
\ep
\bpf
Note that by the definition of the central character $\tomega_i$ one has $\tomega_i(\mtr S^\ro) =\tomega_i(\mtc F(F^\ro))=\frac{\lag\mtc F(F^\ro), \;x_i\rag}{\fp(x_i)}.$
Therefore we have to show that
$\lag\mtc F(F^\ro), x_i\rag=\ro(x_i)\frac{\fpr}{ c_\ro}$ for all $x_i$, 
i.e. 
$
\mtc F(F^\ro)=\frac{\fpr}{ c_\ro}\ro.
$
By applying $\mtc G$ to the above equation, since $\mtc G$ is a bijection and  $\mtc G\circ \mtc F=(-)^{*}$, it is enough to show 
$(F^\ro)^*= \frac{\fpr}{ c_\ro}\mtc G(\ro)$. Thus one has to show $\lag\nu, (F^\ro)^*\rag=\frac{\fpr}{ c_\ro}\lag\nu ,\;\mtc G(\ro)\rag, \;\text{for all}\;\nu \in \wdr.$ By definition of $\mtc G$ one has
\begin{eqnarray*}
\lag\nu, \;\frac{\fpr}{ c_\ro}\mtc G(\ro)\rag&=&\frac{\fpr}{ c_\ro}\lag\ro\star \nu, \;F_0\rag=
\\&=& \frac{\fpr}{ c_\ro}\lag\ro\star \nu, \frac{1}{\fpr}\sumitom d_ix_i\rag
\\&=& \frac{1}{ c_\ro}\sumitom\ro(x_i)\nu(x_i)
\end{eqnarray*}
On the other hand by applying $\id \otimes \rho$ to Equation \eqref{db} we get that $F^\ro=\frac{1}{n_\ro}\big(\sumitom\ro(x_{i^*})x_i\big)$
and therefore,
$$
(F^\ro)^*=\frac{1}{n_\ro}\big(\sumitom \ro(x_{i^*})x_{i^*}\big)=\frac{1}{n_\ro}\big(\sumitom\ro(x_{i})x_{i}\big).
$$
This shows that 
$$
\lag\nu, \; (F^\ro)^*\rag=\frac{1}{n_\ro}\big(\sumitom\ro(x_{i})\nu(x_{i})\big)=\frac{1}{ c_\ro}\lag\nu, \mtc G(\ro)\rag
$$
and the proof is finished.
\epf
\bn{defn}\label{ais} \label{def:isa}
We say that a fusion ring $(R, \mtc B)$ is Isaacs if 
$$
\tomega_i(\mtr S^\rho)\in \mathbb A, 
$$
for all $x_i\in \mtc B$ and all $\ro \in \irr(\rc)$.
\end{defn}
\bt\label{is-abstr}
Let $(R,B)$ be a fusion ring with Isaacs property. If $R$ is Isaacs and the element $x_i\in \mtc B$ is central in $R$(i.e. $x_i\in Z(R)$) then 
$$
\frac{\fpr}{\fp(x_i)}\in \mathbb A.
$$
\et
\bpf
Since $x_i\in Z(R)$ one has that 
$$
x_i=\sum_{\ro\in\irr(\rc)} \al_{i\rho}F^\rho
$$
for some scalars $\aliro \in \comp$. Since $\aliro$ are eigenvalues of $L_{x_i}$ they are algebraic integers. Applying the Fourier transform $\mtc F$ to the above equality  it follows that
$$
 \frac{\fpr}{\fp(x_{i^*})}\widetilde E_{i^*}=\sum_{\ro\in\irr(\rc)} \aliro\mtr S^\rho.
$$
Applying moreover $\tomega_{i^*}$ to this new equality it follows that
$$
\frac{\fpr}{\fp(x_{i^*})}=\sum_{\ro\in\irr(\rc)} \aliro\tomega_i(\mtr S^\rho)\in \mathbb A.
$$
since $\tomega_i(\mtr S^\rho)\in \mathbb A$ by the Isaacs property of $R$.
\epf
\subsection{Isaacs for Grothendieck rings}
Suppose that $R=K_0(\cc)$  is the Grothendieck ring of a pseudo-unitary fusion category $\cc$.

Recall  \cite[Definition 5.1]{eno-nec} that $\cc$ has $s$-Issacs property if for any simple object $X$ of $\cc$ and any $\ro \in \irr(K(\cc))$ one $$
\lam_s(\rho, X)\in \mathbb A.
$$   
Here 
$$\lam_s(\rho, X):=\dim(\cc)^s\dim(Z_\rho)^{1-s}\frac{\ro(X)}{\dim(X)}$$
where $Z_\ro\in \irr(\czcc)$ is the representation corresponding to $\ro$, see \cite[Sect 3]{eno-nec}. The 0-Isaacs property is simply called the Isaacs property of $\cc$ and it was introduced previously in \cite{lpr1, lpr2}.
\bt\label{cateno}
Let $\cc$ be a pseudo-unitary fusion category. With the above notations, for $R=K_0(\cc)$ one has that
\beq\label{3eq}
\lam_0(\rho, X_i)=\tomega_i(\mtr S^\rho).
\eeq
Then $\cc$ has the Isaacs property if and only if $K_0(\cc)$ satisfies the abstract Isaacs property from Definition \ref{ais}. 
\et
\bpf
With the above notations it is easy to see that $\dim(\cc^\rho_1)=\dim(Z_\rho)$. Recall that $\mtr C^\rho:=\mtc F_\alm^{-1}(F^\rho)\in \cecc$ is the {\it matrix-class sum} associated to $\ro$. Note that $\al(\mtr S^\ro)=\al(\mtc F(F^\ro))=\mtc F_\alm^{-1}(F^\ro)=C^\ro.$ 

It follows that $\omega_i(\mtr C^\ro)=\omega_i(\al(\mtr S^\ro))\numeq{\ref{moveomega}}\tomega_i(\mtr S^\ro)$. Therefore for $R=K_0(\cc)$ the Isaacs property can be written as $\omega_i(\mtr C^\ro)\in \mathbb A$. On the other hand note that by its definition one has 
$
\omega_i(\mtr C^\rho)=\frac{\lag\ch_i,\;\mtr C^\ro\rag_z}{\fp(X_i)}\numeq{\ref{mujev}}\frac{\ro(\ch_i)\dim(\cc^\ro_1)}{\fp(X_i)}.
$
This gives that 
\beq\label{roci}
\ro(\ch_i)=\frac{\omega_i(\mtr C^\ro)\fp(X_i)}{\dim(\cc^\ro_1)}
\eeq
Note that by its definition one has
$\lam_0(\rho, X_i)=\dim(\cc^\rho_1)\frac{\ro(X_i)}{\dim(X_i)}.$
Using Equation \eqref{roci} one obtains that 
$\lam_0(\ro, X_i)=\omega_i(\mtr C^\ro)$ since ${\dim(X_i)}={\fp(X_i)}$ in the pseudo-unitary case.
\epf
The following corollary follows from Proposition \ref{is-abstr}.
\bc\label{is-categ}
Let $\cc$ be a pseudo-unitary fusion category with Isaacs property and $X_i$ be a simple object of $\cc$. If $\ch_i\in Z(\cfcc)$ is a central character  then $\frac{\fp(\cc)}{\fp(X_i)}\in \mathbb A$.
\ec

If $\cc$ has a commutative Grothendieck ring then this corollary follows from \cite[Proposition 5.4]{eno-nec}.
\br
It follows from the above proof that if $\cc$ is a pseudo-unitary fusion category then:
$$
\lam_s(\ro, X_i)=(\frac{\dim(\cc)}{\dim(\cc^\ro_1)})^s\omega_i(\mtr C^\ro)=n_\ro^s\omega_i(\mtr C^\ro).
$$
\er
\subsection{Frobenius type results for fusion rings}
For a  fusion ring $(R, B)$ and $\ro \in \irr(\rc)$ recall the matrix-class sum $\mtr S^\rho:=\mtc F(F^\rho)\in \widehta{R_\comp}$.

If $R$ is commutative then $\{F^\rho\}_\ro$ form a linear basis of $\wdr$. Therefore there are some scalars $c^\nu_{\ro, \psi}$ defined by
\beq\label{ssmult}
\mtr S^\rho\star \mtr S^\psi=\sum_{\nu \in \irr(\rc)}c^\nu_{\ro, \psi}\mtrs^\nu.
\eeq 
They are called the {\it structure constants} of $R$.

\bt\label{d-Frobenius}
Let $(R, \mtc B)$ be a commutative fusion ring. Let $d \in \comp$. Consider the following three properties:
\bne
\item[(a)] $dc^\nu_{\ro \psi}\in \mathbb A$ for all $\nu, \ro, \psi$, %\label{f1}
\item[(b)] $\tomega_i(dS^\rho)\in \mathbb A$ for all $i, \rho$, %\label{f2}
\item[(c)] %$R$ satisfies the following Frobenius type property:
$ \frac{d\fpr}{\fp(x_i)}\in \mathbb A, $  for all $i$. %\label{f3}
%for any standard basic element $x_i$ of $\mtc B$.
\ene
Then (a) implies (b), and (b) implies (c).
\et

\bpf
%We will prove that 
%$$
%dc^\nu_{\ro \psi} \in \mathbb A\implies \tomega_i(dC^\rho)\in \mathbb A\implies  \frac{d\fpr}{\fp(x_i)}\in \mathbb A.
%$$
Recall  that $\widetilde E_i={\fp(x_i)}x_i^\circ$ are the central primitive idempotents of $\wdr$ and by Equation \eqref{finvlin} one has $\mtc F(x_i)=\frac{\fpr}{{\fp(x_i)}}\widetilde E_{i^*}$. On the other hand we can write $x_i=\sumroirc\al_{i\rho} F^\ro$
for some complex scalars $\al_{i\rho} \in \mathbb C$. This implies 
$\mtc F(x_i)=\sumroirc \al_{i\rho}\mtrs^\ro$.

Comparing the two equations for $\mtc F(x_i)$ one deduces that
$$\frac{\fpr}{{\fp(x_i)}} \widetilde E_{i^*}=\sum_{\ro \in\irr(\rc)}\al_{i\rho} \mtrs^\ro$$

Note that $\aliro \in \mathbb A$ as eigenvalue of a fusion matrix (having integer entries).  %since $x_ix_j=\sum_l N^l_{ij}x_l$ with $N^l_{ij}\in \mathbb Z_{\geq 0}$.

Recall also that there is an algebra character $\tomega_i:\wdr\ra \comp$, defined by  $\tomega_i(\mu):=\frac{\mu(x_i)}{{\fp(x_i)}}$. It is called the central character assciated to $x_i$. Now we prove (a) $\implies$ (b):

Applying the central character $\tomega_{i^*}$ to the above equation:
\beq\label{one}
\frac{\fpr}{{\fp(x_i)}}=\omega_{i^*}(\frac{\fpr}{{\fp(x_i)}}\widetilde E_{i^*})=\omega_{i^*}(\mtc F(x_i))=\sumroirc\aliro\omega_{i^*}(\mtr S^\ro)
\eeq
Equation \eqref{ssmult} implies that
$$(dS^\ro)(dS^\psi)=\sum_{\nu\in \irr(\rc)} (dc^\nu_{\ro \psi})(dS^\nu)$$ and by a standard argument, see \cite[Theorem 3.4]{is} one has $\omega_i(dS^\nu)\in \mathbb A$.

Finally,  we prove (b) $\implies$ (c): from Equation \eqref{one}  one has
\beqn
\frac{d\fpr}{{\fp(x_i)}}=\sumroirc\aliro\big[\tomega_{i^*}(dS^\ro)\big]\in \mathbb A,
\eeqn
which finishes the proof of the theorem.
\epf

\begin{rem}
For $d=1$, the three properties of Theorem \ref{d-Frobenius} are 
\bne
\item[(a)] algebraic integrality of the structure constants,
\item[(b)] Isaacs property,
\item[(c)] $1$-Frobenius property.
\ene
By \cite{CW6}, the example $\mathcal{Z}(\mathrm{Vec}_{S_3})$ shows that (a) is in fact \emph{strictly} stronger than (b), and by Proposition \ref{prop:EH}, the example $\mathcal{EH}_1$ shows that  (b) is \emph{strictly} stronger than (c). 
\end{rem}

\subsection{About Morita equivalence} \label{sub:MoEq}
This subsection provides a sufficient condition (involving the Morita equivalence) for a property to be true for every spherical fusion category. We deduce that the $1$-Frobenius property holds for every spherical fusion category if and only if it is invariant by Morita equivalence. Idem for the Isaacs property, so that it cannot be invariant by Morita equivalence (by Proposition \ref{prop:EH}).
\bp \label{prop:P}
Let (P) be a property on spherical fusion categories such that:
\begin{itemize}
\item[(1)] it holds for every modular fusion category,
\item[(2)] for every spherical fusion category $\mathcal{C}$, if $\mathcal{C}$ is non-(P) then so is $\mathcal{C} \boxtimes \mathcal{C}^{op}$,
\item[(3)] it is invariant by Morita equivalence.
\end{itemize}
Then (P) holds for every spherical fusion category.
\ep
\bpf
Let $\mathcal{C}$ be a spherical fusion category. By \cite[Theorem 4.24]{MuII}, $\mathcal{Z}(\mathcal{C})$ is Morita equivalent to $\mathcal{C} \boxtimes \mathcal{C}^{op}$. If $\mathcal{C}$ is non-(P) then so is $\mathcal{C} \boxtimes \mathcal{C}^{op}$ by (2), but by (3), (P) is invariant by Morita equivalence, so $\mathcal{Z}(\mathcal{C})$ is also non-(P), but $\mathcal{Z}(\mathcal{C})$ is modular, contradiction with (1).
\epf

\bc \label{cor:FrMo}
The $1$-Frobenius property holds for every spherical fusion category if and only if it is invariant by Morita equivalence.
\ec
\bpf
One way is trivial, and the other way follows from Proposition \ref{prop:P}. More precisely, the $1$-Frobenius property satisfies (1), because it holds more generally for every spherical braided fusion category, see \cite[Corollary 9.3.5]{EGNO15}. About (2), let $\mathcal{C}$ be a spherical fusion category and $X$ a simple object such that $\frac{\dim(\mathcal{C})}{\dim(X)} \not \in \mathbb{A}$, then $X \boxtimes X$ is a simple object of $\mathcal{C} \boxtimes \mathcal{C}^{op}$ and $$\frac{\dim(\mathcal{C} \boxtimes \mathcal{C}^{op})}{\dim(X \boxtimes X)} = \left(\frac{\dim(\mathcal{C})}{\dim(X)} \right)^2 \not \in \mathbb{A}.$$
Finally (3) holds by assumption.
\epf

By \emph{Isaacs terms} below we mean the numbers in Definition \ref{def:isa}.
\bl \label{lem:IsaDe} 
Let $\mathcal{C}$ and $\mathcal{D}$ be two spherical fusion categories. Then the Isaacs terms of $\mathcal{C} \boxtimes \mathcal{D}$ are the products of the Isaacs terms of $\mathcal{C}$ and $\mathcal{D}$. In particular, the Isaacs property is invariant by Deligne tensor product.

\el
\bpf
Let $R$ and $S$ be the Grothendieck rings of $\mathcal{C}$ and $\mathcal{D}$ respectively. Then the Grothendieck ring of $\mathcal{C} \boxtimes \mathcal{D}$ is  $R \otimes S$, but the Fourier transform on $R \otimes S$ is the tensor product of the Fourier transforms on $R$ and $S$. The result follows.
\epf

\bc
The Isaacs property is not invariant by Morita equivalence.
\ec
\bpf
Let us apply Proposition \ref{prop:P}: the Isaacs property satisfies (1) by \cite[Proposition 5.2]{eno-nec}, and it satisfies (2) as for the proof of Corollary \ref{cor:FrMo}, by Lemma \ref{lem:IsaDe}. Finally, if it satisfies (3), i.e. if the Isaacs property is invariant by Morita equivalence, then it is true for every spherical fusion category, contradiction with Proposition \ref{prop:EH}.
\epf

\section{Extended Haagerup fusion categories} \label{sec:EH}
The section provides some data about the Extended Haagerup fusion categories and proves the following result.  
\begin{prop} \label{prop:EH} The Extended Haagerup fusion categories are non-Isaacs.
\end{prop}
The fusion matrices come from \cite{gmpps}, and we deduced the other data using SageMath.
\begin{notations} \label{def:nota} 
Let $p$ be an odd prime, $\zeta_p:= e^{2i\pi/p}$ and $m:= (p-1)/2$. Let
$[a_1,\dots, a_m]_p := -\sum_{k=1}^{m} a_k(\zeta_p^k+\zeta_p^{-k})$.
\end{notations}
\subsection{Some data about $\mathcal{EH}_1$}  \label{data}
Here are some data about the fusion category $\mathcal{EH}_1$. 

\begin{itemize}
\item Fusion matrices $M_1, M_2, \dots, M_6$: 
$$ \begin{smallmatrix}1 & 0 & 0 & 0 & 0 & 0 \\0 & 1 & 0 & 0 & 0 & 0 \\0 & 0 & 1 & 0 & 0 & 0 \\0 & 0 & 0 & 1 & 0 & 0 \\0 & 0 & 0 & 0 & 1 & 0 \\0 & 0 & 0 & 0 & 0 & 1\end{smallmatrix}  , \ \begin{smallmatrix}0 & 1 & 0 & 0 & 0 & 0 \\1 & 1 & 1 & 0 & 0 & 0 \\0 & 1 & 1 & 0 & 1 & 0 \\0 & 0 & 0 & 0 & 1 & 1 \\0 & 0 & 1 & 1 & 1 & 1 \\0 & 0 & 0 & 1 & 1 & 2\end{smallmatrix}  , \ \begin{smallmatrix}0 & 0 & 1 & 0 & 0 & 0 \\0 & 1 & 1 & 0 & 1 & 0 \\1 & 1 & 1 & 1 & 1 & 1 \\0 & 0 & 1 & 1 & 1 & 2 \\0 & 1 & 1 & 1 & 2 & 3 \\0 & 0 & 1 & 2 & 3 & 3\end{smallmatrix}  , \ \begin{smallmatrix}0 & 0 & 0 & 1 & 0 & 0 \\0 & 0 & 0 & 0 & 1 & 1 \\0 & 0 & 1 & 1 & 1 & 2 \\1 & 0 & 1 & 1 & 2 & 2 \\0 & 1 & 1 & 2 & 3 & 3 \\0 & 1 & 2 & 2 & 3 & 4\end{smallmatrix}  , \ \begin{smallmatrix}0 & 0 & 0 & 0 & 1 & 0 \\0 & 0 & 1 & 1 & 1 & 1 \\0 & 1 & 1 & 1 & 2 & 3 \\0 & 1 & 1 & 2 & 3 & 3 \\1 & 1 & 2 & 3 & 4 & 5 \\0 & 1 & 3 & 3 & 5 & 6\end{smallmatrix}  , \ \begin{smallmatrix}0 & 0 & 0 & 0 & 0 & 1 \\0 & 0 & 0 & 1 & 1 & 2 \\0 & 0 & 1 & 2 & 3 & 3 \\0 & 1 & 2 & 2 & 3 & 4 \\0 & 1 & 3 & 3 & 5 & 6 \\1 & 2 & 3 & 4 & 6 & 7\end{smallmatrix} $$
\item $\fp$: $$[170, 120, 120, 295, 170, 295]_{13} \simeq 570.246818815795.$$
\item Type $[d_1, d_2, \dots, d_6]$ where $d_i:=\fp(X_i)$: \begin{align*}
&\left[ 1,[2,1,1,2,2,2]_{13},[3,2,2,4,3,4]_{13},[2,1,1,4,2,4]_{13},[4,3,3,7,4,7]_{13},[4,3,3,8,4,8]_{13}\right]   \\
\simeq &\left[ 1,
 3.37720 2853972,
 7.028296262910,
 8.679389671847,
 13.33048308078,
 15.98157648972 \right].
 \end{align*}
 \item Formal codegrees $[c_1, c_2, \dots, c_6]$: \begin{align*}
&\left[ [170, 120, 120, 295, 170, 295]_{13},[120, 295, 295, 170, 120, 170]_{13},5,5,5,[295, 170, 170, 120, 295, 120]_{13} \right]  \\ 
 \simeq &[570.246818815795,
 11.5441710015915,
 5,
 5,
 5,
 3.20901018261429]
 \end{align*}

\item Character table $[\lambda_{i,j}]_{i,j \in \{1,\dots, 6\}}$:
\begin{align*} &\left[ \begin{matrix}1 & 1 & 1 & 1 & 1 & 1 \\
 [2,1,1,2,2,2]_{13} & [1,2,2,2,1,2]_{13} & 1 & -1\phantom{-} & 0 & [2,2,2,1,2,1]_{13} \\
  [3,2,2,4,3,4]_{13} & [2,4,4,3,2,3]_{13} & -1\phantom{-} & 1 & -1\phantom{-} & [4,3,3,2,4,2]_{13} \\
   [2,1,1,4,2,4]_{13} & [1,4,4,2,1,2]_{13} & 0 & 1 & 1 & [4,2,2,1,4,1]_{13} \\
    [4,3,3,7,4,7]_{13} & [3,7,7,4,3,4]_{13} & -1\phantom{-} & -1\phantom{-} & 1 & [7,4,4,3,7,3]_{13} \\
     [4,3,3,8,4,8]_{13} & [3,8,8,4,3,4]_{13} & 1 & 0 & -1\phantom{-} & [8,4,4,3,8,3]_{13}
\end{matrix} \right]   \\
\simeq &\left[ \begin{matrix}1 & 1 & 1 & 1 & 1 & 1 \\ 3.37720285397 & 2.27389055496 & 1 & -1\phantom{-} & 0 & -0.651093408937 \\ 7.02829626291 & 1.89668770099 & -1\phantom{-} & 1 & -1\phantom{-} & 0.0750160360986 \\ 8.67938967184 & -0.48051515298 & 0 & 1 & 1 & -1.19887451886 \\ 13.3304830807 & 0.142281993045 & -1\phantom{-} & -1\phantom{-} & 1 & 0.527234926170 \\ 15.9815764897 & -1.23492086092 & 1 & 0 & -1\phantom{-} & 0.253344371205\end{matrix} \right]
 \end{align*}
\end{itemize}

Note that all the data can be read in the character table. The type is the first column, the formal codegrees are the squared norm of the columns, where the biggest one is $\fp$.

\subsection{$\mathcal{EH}_1$ is not Isaacs} \label{sub:EH1}
Here we will show that the fusion category $\mathcal{EH}_1$ does not satisfies the Isaacs property. We will use the notations $[a_1,\dots, a_m]_{13}$, $M_i$, $d_i$, $c_i$ and $\lambda_{i,j}$ from \S \ref{data}. Take $i=j=2$, then 
\begin{align*}
\lambda_{i,j} &= [1,2,2,2,1,2]_{13},\\
c_1 &= [170,120,120,295,170,295]_{13},\\
d_i &=[2,1,1,2,2,2]_{13}, \\
c_j &=[120,295,295,170,120,170]_{13}, \\
\frac{\lambda_{i,j}c_1}{d_ic_j} &= [9,\frac{32}{5},\frac{32}{5},\frac{84}{5},9,\frac{84}{5}]_{13} \not \in \mathbb{Z}[\zeta_{13}],
\end{align*}
so the Isaacs property is not satisfied. 

\begin{remark}
Note that the $s$-Isaacs property introduced in \cite{eno-nec} means that for all $i,j$ then
$\frac{\lambda_{i,j}c_1}{d_ic_j^{1-s}}$ is an algebraic integer. Then $0$-Isaacs corresponds to Isaacs. Observe that $\mathcal{EH}_1$ is $s$-Isaacs if and only if $s \ge 1$.
\end{remark}

\subsection{$\mathcal{EH}_i$ are not Isaacs} For $i=2,3,4$, the Grothendieck rings $R_i$ of $\mathcal{EH}_i$ are noncommutative, more precisely $$R_i \otimes_{\mathbb{Z}}\mathbb{C} \simeq \mathbb{C}^{\oplus 4} \oplus M_2(\mathbb{C}).$$
We show that none is Isaacs. For so, we (luckily) only need to consider the characters of degree $1$ (i.e. the central part). Here are their fusion matrices, and \textbf{in the central part}, their character table and formal codegrees. 

\begin{itemize}
\item $\mathcal{EH}_2$
$$\begin{smallmatrix}
\begin{smallmatrix} 1 & 0 & 0 & 0 & 0 & 0 & 0 & 0 \\ 0 & 1 & 0 & 0 & 0 & 0 & 0 & 0 \\ 0 & 0 & 1 & 0 & 0 & 0 & 0 & 0 \\ 0 & 0 & 0 & 1 & 0 & 0 & 0 & 0 \\ 0 & 0 & 0 & 0 & 1 & 0 & 0 & 0 \\ 0 & 0 & 0 & 0 & 0 & 1 & 0 & 0 \\ 0 & 0 & 0 & 0 & 0 & 0 & 1 & 0 \\ 0 & 0 & 0 & 0 & 0 & 0 & 0 & 1\end{smallmatrix}, \ 
\begin{smallmatrix} 0 & 1 & 0 & 0 & 0 & 0 & 0 & 0 \\ 1 & 1 & 0 & 0 & 1 & 0 & 0 & 0 \\ 0 & 0 & 0 & 0 & 0 & 0 & 1 & 0 \\ 0 & 0 & 0 & 0 & 0 & 1 & 0 & 0 \\ 0 & 1 & 0 & 0 & 1 & 0 & 0 & 1 \\ 0 & 0 & 0 & 1 & 0 & 1 & 1 & 1 \\ 0 & 0 & 1 & 0 & 0 & 1 & 1 & 1 \\ 0 & 0 & 0 & 0 & 1 & 1 & 1 & 1\end{smallmatrix}, \  
\begin{smallmatrix} 0 & 0 & 1 & 0 & 0 & 0 & 0 & 0 \\ 0 & 0 & 0 & 0 & 0 & 1 & 0 & 0 \\ 0 & 0 & 0 & 0 & 0 & 0 & 0 & 1 \\ 1 & 0 & 0 & 0 & 0 & 1 & 0 & 0 \\ 0 & 0 & 0 & 0 & 0 & 0 & 1 & 1 \\ 0 & 0 & 0 & 0 & 1 & 1 & 1 & 1 \\ 0 & 1 & 1 & 0 & 0 & 1 & 1 & 1 \\ 0 & 0 & 0 & 1 & 1 & 1 & 1 & 1\end{smallmatrix}, \  
\begin{smallmatrix} 0 & 0 & 0 & 1 & 0 & 0 & 0 & 0 \\ 0 & 0 & 0 & 0 & 0 & 0 & 1 & 0 \\ 1 & 0 & 0 & 0 & 0 & 0 & 1 & 0 \\ 0 & 0 & 0 & 0 & 0 & 0 & 0 & 1 \\ 0 & 0 & 0 & 0 & 0 & 1 & 0 & 1 \\ 0 & 1 & 0 & 1 & 0 & 1 & 1 & 1 \\ 0 & 0 & 0 & 0 & 1 & 1 & 1 & 1 \\ 0 & 0 & 1 & 0 & 1 & 1 & 1 & 1\end{smallmatrix}, \  
\begin{smallmatrix} 0 & 0 & 0 & 0 & 1 & 0 & 0 & 0 \\ 0 & 1 & 0 & 0 & 1 & 0 & 0 & 1 \\ 0 & 0 & 0 & 0 & 0 & 1 & 0 & 1 \\ 0 & 0 & 0 & 0 & 0 & 0 & 1 & 1 \\ 1 & 1 & 0 & 0 & 1 & 1 & 1 & 1 \\ 0 & 0 & 1 & 0 & 1 & 2 & 2 & 2 \\ 0 & 0 & 0 & 1 & 1 & 2 & 2 & 2 \\ 0 & 1 & 1 & 1 & 1 & 2 & 2 & 2\end{smallmatrix}, \  
\begin{smallmatrix} 0 & 0 & 0 & 0 & 0 & 1 & 0 & 0 \\ 0 & 0 & 1 & 0 & 0 & 1 & 1 & 1 \\ 0 & 1 & 1 & 0 & 0 & 1 & 1 & 1 \\ 0 & 0 & 0 & 0 & 1 & 1 & 1 & 1 \\ 0 & 0 & 0 & 1 & 1 & 2 & 2 & 2 \\ 1 & 1 & 1 & 1 & 2 & 4 & 3 & 3 \\ 0 & 1 & 1 & 1 & 2 & 3 & 3 & 4 \\ 0 & 1 & 1 & 1 & 2 & 3 & 4 & 4\end{smallmatrix}, \  
\begin{smallmatrix} 0 & 0 & 0 & 0 & 0 & 0 & 1 & 0 \\ 0 & 0 & 0 & 1 & 0 & 1 & 1 & 1 \\ 0 & 0 & 0 & 0 & 1 & 1 & 1 & 1 \\ 0 & 1 & 0 & 1 & 0 & 1 & 1 & 1 \\ 0 & 0 & 1 & 0 & 1 & 2 & 2 & 2 \\ 0 & 1 & 1 & 1 & 2 & 3 & 3 & 4 \\ 1 & 1 & 1 & 1 & 2 & 3 & 4 & 3 \\ 0 & 1 & 1 & 1 & 2 & 4 & 3 & 4\end{smallmatrix}, \  
\begin{smallmatrix} 0 & 0 & 0 & 0 & 0 & 0 & 0 & 1 \\ 0 & 0 & 0 & 0 & 1 & 1 & 1 & 1 \\ 0 & 0 & 0 & 1 & 1 & 1 & 1 & 1 \\ 0 & 0 & 1 & 0 & 1 & 1 & 1 & 1 \\ 0 & 1 & 1 & 1 & 1 & 2 & 2 & 2 \\ 0 & 1 & 1 & 1 & 2 & 3 & 4 & 4 \\ 0 & 1 & 1 & 1 & 2 & 4 & 3 & 4 \\ 1 & 1 & 1 & 1 & 2 & 4 & 4 & 4\end{smallmatrix} 
\end{smallmatrix}$$  

$$\left[ \begin{matrix}1 & 1 & 1 & 1 \\ 
 3.37720285397296 & 2.27389055496422 & 0 & -0.651093408937175 \\
 3.65109340893718 & -0.377202853972958 & -1 & 0.726109445035790 \\
 3.65109340893718 & -0.377202853972958 & -1 & 0.726109445035790 \\
 7.02829626291013 & 1.89668770099126 & -1 & 0.0750160360986070 \\
 12.3304830807845 & -0.857718006954660 & 0 & -0.472765073829828 \\
 12.3304830807845 & -0.857718006954660 & 0 & -0.472765073829828 \\
 13.3304830807845 & 0.142281993045350 & 1 & 0.527234926170180 \\
 \hline
 570.246818815795 & 11.5441710015912 & 5 & 3.20901018261404
 \end{matrix} \right]$$
\noindent The Isaacs property fails by considering the entry $\simeq 2.27389055496422$, as for $\mathcal{EH}_1$, we get the Isaacs term $$\frac{[1,2,2,2,1,2]_{13} [170,120,120,295,170,295]_{13}}{[2,1,1,2,2,2]_{13} [120,295,295,170,120,170]_{13}} = [9,\frac{32}{5},\frac{32}{5},\frac{84}{5},9,\frac{84}{5}]_{13} \not \in \mathbb{Z}[\zeta_{13}].$$

\item $\mathcal{EH}_3$
$$\begin{smallmatrix}
\begin{smallmatrix} 1 & 0 & 0 & 0 & 0 & 0 & 0 & 0 \\ 0 & 1 & 0 & 0 & 0 & 0 & 0 & 0 \\ 0 & 0 & 1 & 0 & 0 & 0 & 0 & 0 \\ 0 & 0 & 0 & 1 & 0 & 0 & 0 & 0 \\ 0 & 0 & 0 & 0 & 1 & 0 & 0 & 0 \\ 0 & 0 & 0 & 0 & 0 & 1 & 0 & 0 \\ 0 & 0 & 0 & 0 & 0 & 0 & 1 & 0 \\ 0 & 0 & 0 & 0 & 0 & 0 & 0 & 1\end{smallmatrix}, \  
\begin{smallmatrix} 0 & 1 & 0 & 0 & 0 & 0 & 0 & 0 \\ 0 & 0 & 0 & 0 & 0 & 1 & 0 & 0 \\ 1 & 0 & 0 & 0 & 1 & 0 & 0 & 0 \\ 0 & 1 & 0 & 0 & 0 & 0 & 1 & 0 \\ 0 & 0 & 0 & 0 & 0 & 0 & 0 & 1 \\ 0 & 0 & 1 & 0 & 0 & 0 & 0 & 1 \\ 0 & 0 & 0 & 0 & 1 & 0 & 1 & 1 \\ 0 & 0 & 0 & 1 & 0 & 1 & 1 & 1\end{smallmatrix}, \  
\begin{smallmatrix} 0 & 0 & 1 & 0 & 0 & 0 & 0 & 0 \\ 1 & 0 & 0 & 1 & 0 & 0 & 0 & 0 \\ 0 & 0 & 0 & 0 & 0 & 1 & 0 & 0 \\ 0 & 0 & 0 & 0 & 0 & 0 & 0 & 1 \\ 0 & 0 & 1 & 0 & 0 & 0 & 1 & 0 \\ 0 & 1 & 0 & 0 & 0 & 0 & 0 & 1 \\ 0 & 0 & 0 & 1 & 0 & 0 & 1 & 1 \\ 0 & 0 & 0 & 0 & 1 & 1 & 1 & 1\end{smallmatrix}, \  
\begin{smallmatrix} 0 & 0 & 0 & 1 & 0 & 0 & 0 & 0 \\ 0 & 0 & 0 & 0 & 0 & 0 & 0 & 1 \\ 0 & 0 & 1 & 0 & 0 & 0 & 1 & 0 \\ 1 & 0 & 0 & 1 & 0 & 0 & 1 & 1 \\ 0 & 0 & 0 & 0 & 0 & 1 & 1 & 1 \\ 0 & 0 & 0 & 0 & 1 & 1 & 1 & 1 \\ 0 & 0 & 1 & 1 & 1 & 1 & 2 & 2 \\ 0 & 1 & 0 & 1 & 1 & 1 & 2 & 3\end{smallmatrix}, \  
\begin{smallmatrix} 0 & 0 & 0 & 0 & 1 & 0 & 0 & 0 \\ 0 & 1 & 0 & 0 & 0 & 0 & 1 & 0 \\ 0 & 0 & 0 & 0 & 0 & 0 & 0 & 1 \\ 0 & 0 & 0 & 0 & 0 & 1 & 1 & 1 \\ 1 & 0 & 0 & 0 & 1 & 0 & 1 & 1 \\ 0 & 0 & 0 & 1 & 0 & 1 & 1 & 1 \\ 0 & 1 & 0 & 1 & 1 & 1 & 2 & 2 \\ 0 & 0 & 1 & 1 & 1 & 1 & 2 & 3\end{smallmatrix}, \  
\begin{smallmatrix} 0 & 0 & 0 & 0 & 0 & 1 & 0 & 0 \\ 0 & 0 & 1 & 0 & 0 & 0 & 0 & 1 \\ 0 & 1 & 0 & 0 & 0 & 0 & 0 & 1 \\ 0 & 0 & 0 & 0 & 1 & 1 & 1 & 1 \\ 0 & 0 & 0 & 1 & 0 & 1 & 1 & 1 \\ 1 & 0 & 0 & 1 & 1 & 1 & 1 & 1 \\ 0 & 0 & 0 & 1 & 1 & 1 & 2 & 3 \\ 0 & 1 & 1 & 1 & 1 & 1 & 3 & 3\end{smallmatrix}, \  
\begin{smallmatrix} 0 & 0 & 0 & 0 & 0 & 0 & 1 & 0 \\ 0 & 0 & 0 & 1 & 0 & 0 & 1 & 1 \\ 0 & 0 & 0 & 0 & 1 & 0 & 1 & 1 \\ 0 & 1 & 0 & 1 & 1 & 1 & 2 & 2 \\ 0 & 0 & 1 & 1 & 1 & 1 & 2 & 2 \\ 0 & 0 & 0 & 1 & 1 & 1 & 2 & 3 \\ 1 & 1 & 1 & 2 & 2 & 2 & 4 & 5 \\ 0 & 1 & 1 & 2 & 2 & 3 & 5 & 6\end{smallmatrix}, \  
\begin{smallmatrix} 0 & 0 & 0 & 0 & 0 & 0 & 0 & 1 \\ 0 & 0 & 0 & 0 & 1 & 1 & 1 & 1 \\ 0 & 0 & 0 & 1 & 0 & 1 & 1 & 1 \\ 0 & 0 & 1 & 1 & 1 & 1 & 2 & 3 \\ 0 & 1 & 0 & 1 & 1 & 1 & 2 & 3 \\ 0 & 1 & 1 & 1 & 1 & 1 & 3 & 3 \\ 0 & 1 & 1 & 2 & 2 & 3 & 5 & 6 \\ 1 & 1 & 1 & 3 & 3 & 3 & 6 & 7\end{smallmatrix} 
\end{smallmatrix}$$

$$\left[ \begin{matrix}1 & 1 & 1 & 1 \\ 
 2.65109340893718 & -1.37720285397296 & 1 & -0.273890554964218 \\
 2.65109340893718 & -1.37720285397296 & 1 & -0.273890554964218 \\
 6.02829626291014 & 0.896687700991260 & 0 & -0.924983963901393 \\
 6.02829626291014 & 0.896687700991260 & 0 & -0.924983963901393 \\
 7.02829626291014 & 1.89668770099126 & 1 & 0.0750160360986070 \\
 13.3304830807845 & 0.142281993045350 & -1 & 0.527234926170180 \\
 15.9815764897217 & -1.23492086092762 & 0 & 0.253344371205960 \\
 \hline
 570.246818815795 & 11.5441710015912 & 5 & 3.20901018261404
 \end{matrix} \right]$$
The Isaacs property fails by considering the entry $\simeq -1.23492086092761$, we get the Isaacs term
$$\frac{[3,8,8,4,3,4]_{13} [170,120,120,295,170,295]_{13}}{[4,3,3,8,4,8]_{13} [120, 295, 295, 170, 120, 170]_{13}} = -[\frac{7}{5}, \frac{4}{5}, \frac{4}{5}, 2, \frac{7}{5}, 2]_{13}  \not \in \mathbb{Z}[\zeta_{13}].$$

\item $\mathcal{EH}_4$
$$\begin{smallmatrix}
\begin{smallmatrix} 1 & 0 & 0 & 0 & 0 & 0 & 0 & 0 \\ 0 & 1 & 0 & 0 & 0 & 0 & 0 & 0 \\ 0 & 0 & 1 & 0 & 0 & 0 & 0 & 0 \\ 0 & 0 & 0 & 1 & 0 & 0 & 0 & 0 \\ 0 & 0 & 0 & 0 & 1 & 0 & 0 & 0 \\ 0 & 0 & 0 & 0 & 0 & 1 & 0 & 0 \\ 0 & 0 & 0 & 0 & 0 & 0 & 1 & 0 \\ 0 & 0 & 0 & 0 & 0 & 0 & 0 & 1\end{smallmatrix}, \  
\begin{smallmatrix} 0 & 1 & 0 & 0 & 0 & 0 & 0 & 0 \\ 0 & 0 & 0 & 1 & 0 & 1 & 1 & 1 \\ 1 & 0 & 0 & 1 & 1 & 0 & 1 & 1 \\ 0 & 1 & 1 & 0 & 1 & 0 & 1 & 1 \\ 0 & 1 & 0 & 1 & 1 & 1 & 1 & 1 \\ 0 & 1 & 1 & 1 & 1 & 1 & 1 & 1 \\ 0 & 0 & 1 & 0 & 1 & 1 & 1 & 2 \\ 0 & 1 & 1 & 1 & 1 & 2 & 1 & 2\end{smallmatrix}, \  
\begin{smallmatrix} 0 & 0 & 1 & 0 & 0 & 0 & 0 & 0 \\ 1 & 0 & 0 & 1 & 1 & 1 & 0 & 1 \\ 0 & 0 & 0 & 1 & 0 & 1 & 1 & 1 \\ 0 & 1 & 1 & 0 & 1 & 1 & 0 & 1 \\ 0 & 0 & 1 & 1 & 1 & 1 & 1 & 1 \\ 0 & 1 & 0 & 0 & 1 & 1 & 1 & 2 \\ 0 & 1 & 1 & 1 & 1 & 1 & 1 & 1 \\ 0 & 1 & 1 & 1 & 1 & 1 & 2 & 2\end{smallmatrix}, \  
\begin{smallmatrix} 0 & 0 & 0 & 1 & 0 & 0 & 0 & 0 \\ 0 & 1 & 1 & 0 & 1 & 1 & 0 & 1 \\ 0 & 1 & 1 & 0 & 1 & 0 & 1 & 1 \\ 1 & 0 & 0 & 1 & 1 & 1 & 1 & 1 \\ 0 & 1 & 1 & 1 & 1 & 1 & 1 & 1 \\ 0 & 1 & 0 & 1 & 1 & 1 & 1 & 2 \\ 0 & 0 & 1 & 1 & 1 & 1 & 1 & 2 \\ 0 & 1 & 1 & 1 & 1 & 2 & 2 & 2\end{smallmatrix}, \  
\begin{smallmatrix} 0 & 0 & 0 & 0 & 1 & 0 & 0 & 0 \\ 0 & 1 & 0 & 1 & 1 & 1 & 1 & 1 \\ 0 & 0 & 1 & 1 & 1 & 1 & 1 & 1 \\ 0 & 1 & 1 & 1 & 1 & 1 & 1 & 1 \\ 1 & 1 & 1 & 1 & 1 & 1 & 1 & 2 \\ 0 & 1 & 1 & 1 & 1 & 2 & 1 & 2 \\ 0 & 1 & 1 & 1 & 1 & 1 & 2 & 2 \\ 0 & 1 & 1 & 1 & 2 & 2 & 2 & 3\end{smallmatrix}, \  
\begin{smallmatrix} 0 & 0 & 0 & 0 & 0 & 1 & 0 & 0 \\ 0 & 0 & 1 & 0 & 1 & 1 & 1 & 2 \\ 0 & 1 & 1 & 1 & 1 & 1 & 1 & 1 \\ 0 & 0 & 1 & 1 & 1 & 1 & 1 & 2 \\ 0 & 1 & 1 & 1 & 1 & 2 & 1 & 2 \\ 1 & 1 & 1 & 1 & 2 & 1 & 2 & 2 \\ 0 & 1 & 1 & 1 & 1 & 2 & 2 & 2 \\ 0 & 2 & 1 & 2 & 2 & 2 & 2 & 3\end{smallmatrix}, \  
\begin{smallmatrix} 0 & 0 & 0 & 0 & 0 & 0 & 1 & 0 \\ 0 & 1 & 1 & 1 & 1 & 1 & 1 & 1 \\ 0 & 1 & 0 & 0 & 1 & 1 & 1 & 2 \\ 0 & 1 & 0 & 1 & 1 & 1 & 1 & 2 \\ 0 & 1 & 1 & 1 & 1 & 1 & 2 & 2 \\ 0 & 1 & 1 & 1 & 1 & 2 & 2 & 2 \\ 1 & 1 & 1 & 1 & 2 & 2 & 1 & 2 \\ 0 & 1 & 2 & 2 & 2 & 2 & 2 & 3\end{smallmatrix}, \  
\begin{smallmatrix} 0 & 0 & 0 & 0 & 0 & 0 & 0 & 1 \\ 0 & 1 & 1 & 1 & 1 & 1 & 2 & 2 \\ 0 & 1 & 1 & 1 & 1 & 2 & 1 & 2 \\ 0 & 1 & 1 & 1 & 1 & 2 & 2 & 2 \\ 0 & 1 & 1 & 1 & 2 & 2 & 2 & 3 \\ 0 & 1 & 2 & 2 & 2 & 2 & 2 & 3 \\ 0 & 2 & 1 & 2 & 2 & 2 & 2 & 3 \\ 1 & 2 & 2 & 2 & 3 & 3 & 3 & 4\end{smallmatrix} 
\end{smallmatrix}$$ 
\end{itemize}

$$\left[ \begin{matrix}1 & 1 & 1 & 1 \\ 
 6.30218681787435 & -1.75440570794592 & 0 & 0.452218890071565 \\
 6.30218681787435 & -1.75440570794592 & 0 & 0.452218890071565 \\
 7.02829626291013 & 1.89668770099126 & -1 & 0.0750160360986070 \\
 8.67938967184731 & -0.480515152981698 & 0 & -1.19887451886561 \\
 9.67938967184731 & 0.519484847018302 & 1 & -0.198874518865610 \\
 9.67938967184731 & 0.519484847018302 & 1 & -0.198874518865610 \\
 13.3304830807845 & 0.142281993045350 & -1 & 0.527234926170180 \\
 \hline
 570.246818815795 & 11.5441710015912 & 5 & 3.20901018261404
 \end{matrix} \right]$$
 The Isaacs property fails by considering the entry $\simeq -0.480515152981694$, we get the Isaacs term
$$\frac{[1,4,4,2,1,2]_{13} [170,120,120,295,170,295]_{13}}{[2,1,1,4,2,4]_{13} [120, 295, 295, 170, 120, 170]_{13}} = [\frac{19}{5},2,2,\frac{2}{5},\frac{19}{5},\frac{2}{5}]_{13} \not \in \mathbb{Z}[\zeta_{13}].$$
\vskip .25cm
\subsection*{Availability of data and materials} 
The datasets generated during and/or analyzed during the current study are available from the second author on reasonable request.
\vskip .25cm
\noindent {\bf \large Declaration}
\vskip .25cm
\noindent {\bf Conflict of interests:}  The authors declare that they have no conflict of interest.
\bibliographystyle{plain}
\bibliography{24nov}

\begin{thebibliography}{10}

\bibitem{ccc-march}
S.~Burciu.
\newblock \textnormal{Conjugacy classes and centralizers for pivotal fusion
  categories}.
\newblock {\em Monatshefte f\"ur Mathematik,}, 193(2):13--46, 2020.

\bibitem{b-blms}
S.~Burciu.
\newblock \textnormal{Structure constants for premodular categories}.
\newblock {\em Bull. Lond. Math. Soc.}, 53(3):777--791, 2021.

\bibitem{b-etale}
S.~Burciu.
\newblock \textnormal{Subalgebras of etale algebras and fusion subcategories}.
\newblock {\em arXiv:2105.11202}, 2021.

\bibitem{b-jpaa}
S.~Burciu.
\newblock \textnormal{On the Galois symmetries for the character table of an
  integral fusion category}.
\newblock {\em J. Algebra Appl.}, 22(1):Paper No. 2350026, 2023.

\bibitem{CW6}
M.~Cohen and S.~Westreich.
\newblock \textnormal{Conjugacy Classes, Class Sums and Character Tables for
  Hopf Algebras}.
\newblock {\em Comm. Algebra}, 39(12):4618--4633, 200.

\bibitem{EGNO15}
P.~Etingof, S.~Gelaki, D.~Nikshych, and V.~Ostrik.
\newblock {\em \textnormal{Tensor categories}}, volume 205.
\newblock Mathematical Surveys and Monographs, American Mathematical Society,
  Providence, RI, 2015.

\bibitem{eno-weakly}
P.~Etingof, D.~Nikshych, and V.~Ostrik.
\newblock \textnormal{Weakly group-theoretical and solvable fusion categories}.
\newblock {\em Adv. Math.}, 226(1):176--205, 2011.

\bibitem{eno-nec}
P.~Etingof, D.~Nikshych, and V.~Ostrik.
\newblock \textnormal{On a necessary condition for unitary categorification of
  fusion rings}.
\newblock {\em arXiv:2102.13239}, 2021.

\bibitem{hdk}
D.~K. Harrison.
\newblock \textnormal{Double coset and orbit spaces}.
\newblock {\em Pacific J. of Math.}, 80(2):451--491, 1979.

\bibitem{hlpw}
L.~Huang, Z.~Liu, S.~Palcoux, and J.~Wu.
\newblock \textnormal{Complete Positivity of Comultiplication and Primary
  Criteria for Unitary Categorification}.
\newblock {\em Int. Math. Res. Not.}, rnad214, 2023.

\bibitem{is}
M.~Isaacs.
\newblock {\em Character theory of finite Groups}.
\newblock Academis Press, New York, San Francisco, London, 1976.

\bibitem{lpr1}
Z.~Liu, S.~Palcoux, and Y.~Ren.
\newblock \textnormal{Classification of Grothendieck rings of complex fusion
  categories of multiplicity one up to rank six}.
\newblock {\em Lett Math Phys}, 112(3):Paper No. 54, 2022.

\bibitem{lpr3}
Z.~Liu, S.~Palcoux, and Y.~Ren.
\newblock \textnormal{Triangular prism equations and categorification}.
\newblock {\em arXiv:2203.06522}, 2022.

\bibitem{lpr2}
Z.~Liu, S.~Palcoux, and Y.~Ren.
\newblock \textnormal{Interpolated family of non-group-like simple integral
  fusion rings of Lie type}.
\newblock {\em Internat. J. Math.}, 34(6):Paper No. 2350030, 2023.

\bibitem{lpw}
Z.~Liu, S.~Palcoux, and J.~Wu.
\newblock \textnormal{Fusion Bialgebras and Fourier Analysis}.
\newblock {\em Adv. Math.}, 390:Paper No. 107905, 2021.

\bibitem{MuII}
M.~M\"{u}ger.
\newblock From subfactors to categories and topology. {II}. {T}he quantum
  double of tensor categories and subfactors.
\newblock {\em J. Pure Appl. Algebra}, 180(1-2):159--219, 2003.

\bibitem{o-fc}
V.~Ostrik.
\newblock \textnormal{On formal codegrees of fusion categories,}.
\newblock {\em Math. Res. Lett.}, 16(5):895--901, 2009.

\bibitem{gmpps}
D.~Penneys E. Peters N.~Snyder P.~Grossman, S.~Morrison.
\newblock \textnormal{The Extended Haagerup fusion categories}.
\newblock {\em Ann. Sci. Éc. Norm. Supér.}, (4) 56(2):589--664, 2023.

\bibitem{MoWa}
K.~Walker S.~Morrison.
\newblock \textnormal{The center of the extended Haagerup subfactor has 22
  simple objects}.
\newblock {\em Internat. J. Math.}, 28(1):1750009, 11, 2017.

\bibitem{scalg}
K.~Shimizu.
\newblock \textnormal{The monoidal center and the character algebra}.
\newblock {\em J. Pure Appl. Alg.}, 221(9):2338--2371, 2017.

\end{thebibliography}
\ed